\documentclass[a4paper, 11pt]{article}
\usepackage{amssymb,latexsym,amsmath,amsthm,amsfonts}
\usepackage{tikz}
\usepackage{paralist}	 
\usepackage{cite}
\usepackage{multicol}
\usepackage{caption}
\usepackage{subfigure}
\usepackage[latin1]{inputenc}
\usepackage{graphicx}
\usepackage{epsf}
\usepackage{fancyhdr}
\usepackage{pstricks}
\usepackage{algorithm}
\usepackage{algorithmic}
\usepackage{lscape}
\newpsobject{showgrid}{psgrid}{subgriddiv=1, griddots=10,gridlabels=6pt}
\usepackage{authblk}
\usepackage[capitalise]{cleveref}

\title{Packing unequal rectangles and squares in a fixed size circular container using formulation space search} 
\author[1]{C.O. López\thanks{claudia.lopez@ciencias.unam.mx}}
\author[2]{J.E. Beasley\thanks{john.beasley@brunel.ac.uk; john.beasley@jbconsultants.biz}}

\affil[1]{Faculty of Sciences, National Autonomous University of Mexico, Mexico City, Mexico}
\affil[2]{Mathematical Sciences, Brunel University, Uxbridge UB8 3PH, UK and JB Consultants, Morden, UK }

\begin{document}
\date{October 2017, Revised February 2018} 
\par\medskip

\maketitle

\begin{abstract}

In this paper we formulate the problem of packing unequal rectangles/squares into a fixed size circular container as a 
mixed-integer nonlinear program. Here we pack rectangles so as to maximise some objective (e.g.~maximise the number of rectangles packed or maximise the total area of the rectangles packed).
 We show how we can eliminate a nonlinear maximisation term that arises in one of the constraints in our formulation. We 
indicate the amendments that can be made to the formulation
for the special case where we are maximising the number of squares packed.
A formulation space search heuristic is presented and
computational results given for publicly available test problems involving up to 30 rectangles/squares. Our heuristic deals with the case where the  rectangles are of fixed orientation (so cannot be rotated) and with the case where the rectangles can be rotated through ninety degrees.

\end{abstract}


\section{Introduction} 

In this paper we consider the problem of packing non-identical rectangles (i.e.~rectangles of different sizes) into a fixed size circular container. Since the circular container may not be large enough to accommodate all of the rectangles available to be packed there exists an element of choice in the problem. In other words we have to decide which of the rectangles will be packed, and moreover for those that are packed their positions within the container. The packing should respect the obvious constraints, namely that the packed rectangles do not overlap with each other  and that each packed rectangle is entirely within the container. This packing should be such so as to maximise an appropriate objective 
(e.g.~maximise the number of rectangles packed or maximise the total area of the rectangles packed).

To illustrate the problem suppose we have ten rectangles with sizes as shown in Table~\ref{table1} to be packed into a fixed sized circular container. The rectangles shown in  Table~\ref{table1} have been ordered into ascending area order.

\begin{table}[!htb]
\centering 

\begin{tabular}{ccc}
\hline
	Rectangle	&	Length	&	Width 	\\   
 \hline
1 & 1.10	& 1.61 	\\  
2 & 2.20	& 1.08	\\  
3 & 1.68	& 1.46	\\  
4 & 1.82	& 2.61	\\  
5 & 2.70	& 2.57	\\  
6 & 3.21	& 2.21	\\  
7 & 2.99	& 3.51	\\  
8 & 3.68	& 3.42	\\  
9 & 4.62	& 3.36	\\  
10 & 3.79	& 4.79	\\  
	\hline	 
											 
\end{tabular}
\caption{Rectangle packing example, circular container radius 4.18}\label{table1}
\end{table}

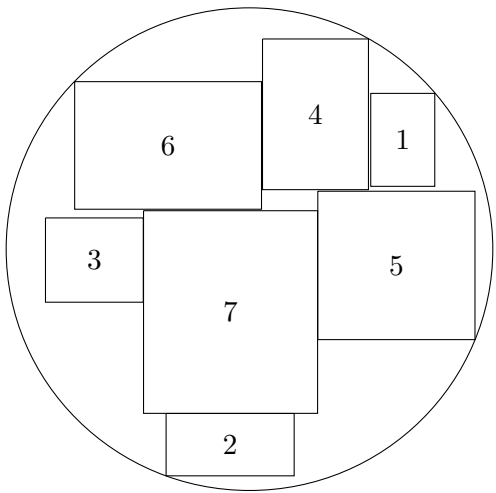
\begin{figure}[!htb]
\centering
\begin{tikzpicture}[scale=0.8];
\draw (0,0) circle [radius=        4.0000000000cm];
\draw(        1.9932421266 ,         2.5819847248) -- 
(        3.0458737055 ,         2.5819847248) -- 
(        3.0458737055 ,         1.0413148683) -- 
(        1.9932421266 ,         1.0413148683) -- 
(        1.9932421266 ,         2.5819847248) ;
 \node at (        2.5195579161 ,         1.8116497966) {     1}; 
\draw(       -1.3708925903 ,        -2.7238255537) -- 
(         .7343705676 ,        -2.7238255537) -- 
(         .7343705676 ,        -3.7573183766) -- 
(       -1.3708925903 ,        -3.7573183766) -- 
(       -1.3708925903 ,        -2.7238255537) ;
 \node at (        -.3182610114 ,        -3.2405719652) {     2}; 
\draw(       -3.3536964542 ,          .5177975242) -- 
(       -1.7460409518 ,          .5177975242) -- 
(       -1.7460409518 ,         -.8793316624) -- 
(       -3.3536964542 ,         -.8793316624) -- 
(       -3.3536964542 ,          .5177975242) ;
 \node at (       -2.5498687030 ,         -.1807670691) {     3}; 
\draw(         .2142232266 ,         3.4840098845) -- 
(        1.9558500209 ,         3.4840098845) -- 
(        1.9558500209 ,          .9864022290) -- 
(         .2142232266 ,          .9864022290) -- 
(         .2142232266 ,         3.4840098845) ;
 \node at (        1.0850366237 ,         2.2352060568) {     4}; 
\draw(        1.1232391583 ,          .9598781146) -- 
(        3.7069712157 ,          .9598781146) -- 
(        3.7069712157 ,        -1.4994520290) -- 
(        1.1232391583 ,        -1.4994520290) -- 
(        1.1232391583 ,          .9598781146) ;
 \node at (        2.4151051870 ,         -.2697869572) {     5}; 
\draw(       -2.8739925679 ,         2.7765491539) -- 
(         .1977777670 ,         2.7765491539) -- 
(         .1977777670 ,          .6617166180) -- 
(       -2.8739925679 ,          .6617166180) -- 
(       -2.8739925679 ,         2.7765491539) ;
 \node at (       -1.3381074004 ,         1.7191328859) {     6}; 
\draw(       -1.7417140510 ,          .6371167752) -- 
(        1.1195299681 ,          .6371167752) -- 
(        1.1195299681 ,        -2.7217348994) -- 
(       -1.7417140510 ,        -2.7217348994) -- 
(       -1.7417140510 ,          .6371167752) ;
 \node at (        -.3110920414 ,        -1.0423090621) {     7}; 
\end{tikzpicture}
\caption{Maximise the number of rectangles packed, no rotation, solution value 7}
\label{fig1}
\end{figure}

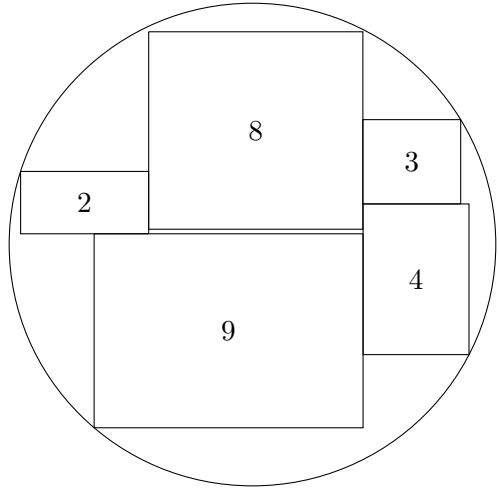
\begin{figure}[!htb]
\centering
\begin{tikzpicture}[scale=0.8];
\draw (0,0) circle [radius=        4.0000000000cm];
\draw(       -3.8119666423 ,         1.2119448296) -- 
(       -1.7067034844 ,         1.2119448296) -- 
(       -1.7067034844 ,          .1784520066) -- 
(       -3.8119666423 ,          .1784520066) -- 
(       -3.8119666423 ,         1.2119448296) ;
 \node at (       -2.7593350633 ,          .6951984181) {     2}; 
\draw(        1.8149071928 ,         2.0702592042) -- 
(        3.4225626952 ,         2.0702592042) -- 
(        3.4225626952 ,          .6731300176) -- 
(        1.8149071928 ,          .6731300176) -- 
(        1.8149071928 ,         2.0702592042) ;
 \node at (        2.6187349440 ,         1.3716946109) {     3}; 
\draw(        1.8179796510 ,          .6730196307) -- 
(        3.5596064452 ,          .6730196307) -- 
(        3.5596064452 ,        -1.8245880248) -- 
(        1.8179796510 ,        -1.8245880248) -- 
(        1.8179796510 ,          .6730196307) ;
 \node at (        2.6887930481 ,         -.5757841970) {     4}; 
\draw(       -1.7066640310 ,         3.5266392263) -- 
(        1.8148670695 ,         3.5266392263) -- 
(        1.8148670695 ,          .2539119536) -- 
(       -1.7066640310 ,          .2539119536) -- 
(       -1.7066640310 ,         3.5266392263) ;
 \node at (         .0541015193 ,         1.8902755899) {     8}; 
\draw(       -2.6031074615 ,          .1782446930) -- 
(        1.8179451701 ,          .1782446930) -- 
(        1.8179451701 ,        -3.0370663117) -- 
(       -2.6031074615 ,        -3.0370663117) -- 
(       -2.6031074615 ,          .1782446930) ;
 \node at (        -.3925811457 ,        -1.4294108093) {     9}; 
\end{tikzpicture}
\caption{Maximise the total area of the rectangles packed, no rotation, solution value 37.6878}
\label{fig2}
\end{figure}

Regarding the rectangles as being of fixed orientation, i.e.~they cannot be rotated, then:
\begin{itemize}

\item If we are wish to maximise the number of rectangles packed  Figure~\ref{fig1} shows the solution as derived by the approach presented in this paper. In that figure we can see that seven of the ten rectangles have been packed, three rectangles are left unpacked.  

\item If we are wish to maximise the total area of the rectangles packed  Figure~\ref{fig2} shows the solution as derived by the approach presented in this paper. In that figure we can see that five of the ten rectangles have been packed.   
\end{itemize}

\begin{figure}[!htb]
\centering
\begin{tikzpicture}[scale=0.8];
\draw (0,0) circle [radius=        4.0000000000cm];
\draw(        1.8712887942 ,        -1.1835039534) -- 
(        2.9239203731 ,        -1.1835039534) -- 
(        2.9239203731 ,        -2.7241738098) -- 
(        1.8712887942 ,        -2.7241738098) -- 
(        1.8712887942 ,        -1.1835039534) ;
 \node at (        2.3976045836 ,        -1.9538388816) {     1     }; 
\draw(       -2.6235913891 ,        -1.6145026574) -- 
(       -1.0159358867 ,        -1.6145026574) -- 
(       -1.0159358867 ,        -3.0116318440) -- 
(       -2.6235913891 ,        -3.0116318440) -- 
(       -2.6235913891 ,        -1.6145026574) ;
 \node at (       -1.8197636379 ,        -2.3130672507) {     3     }; 
\draw(       -3.6064450445 ,          .8593985441) -- 
(       -1.0227129871 ,          .8593985441) -- 
(       -1.0227129871 ,        -1.5999315994) -- 
(       -3.6064450445 ,        -1.5999315994) -- 
(       -3.6064450445 ,          .8593985441) ;
 \node at (       -2.3145790158 ,         -.3702665277) {     5     }; 
\draw(       -2.6445929091 ,         2.9919662931) -- 
(         .4271774259 ,         2.9919662931) -- 
(         .4271774259 ,          .8771337572) -- 
(       -2.6445929091 ,          .8771337572) -- 
(       -2.6445929091 ,         2.9919662931) ;
 \node at (       -1.1087077416 ,         1.9345500251) {     6     }; 
\draw(        -.9992337899 ,          .5050789081) -- 
(        1.8620102292 ,          .5050789081) -- 
(        1.8620102292 ,        -2.8537727665) -- 
(        -.9992337899 ,        -2.8537727665) -- 
(        -.9992337899 ,          .5050789081) ;
 \node at (         .4313882196 ,        -1.1743469292) {     7     }; 
\draw(        2.0597582061 ,          .9321285414) -- 
(        3.0932510290 ,          .9321285414) -- 
(        3.0932510290 ,        -1.1731346165) -- 
(        2.0597582061 ,        -1.1731346165) -- 
(        2.0597582061 ,          .9321285414) ;
 \node at (        2.5765046176 ,         -.1205030376) {     2r    }; 
\draw(         .4470410778 ,         2.6948539495) -- 
(        2.9446487333 ,         2.6948539495) -- 
(        2.9446487333 ,          .9532271552) -- 
(         .4470410778 ,          .9532271552) -- 
(         .4470410778 ,         2.6948539495) ;
 \node at (        1.6958449055 ,         1.8240405524) {     4r    }; 

\end{tikzpicture}
\captionsetup{font=small}
\caption{Maximise the number of rectangles packed, rotation allowed, solution value 7}
\label{fig1a}
\end{figure}
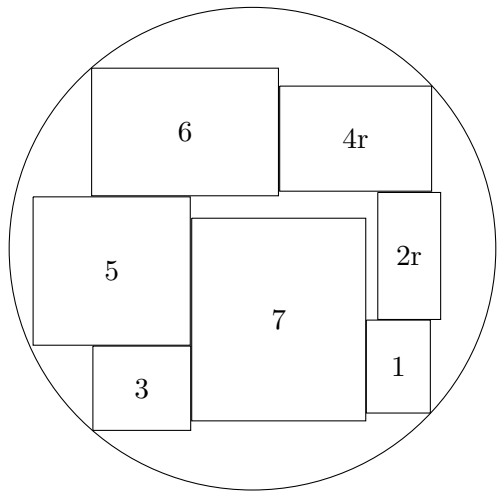

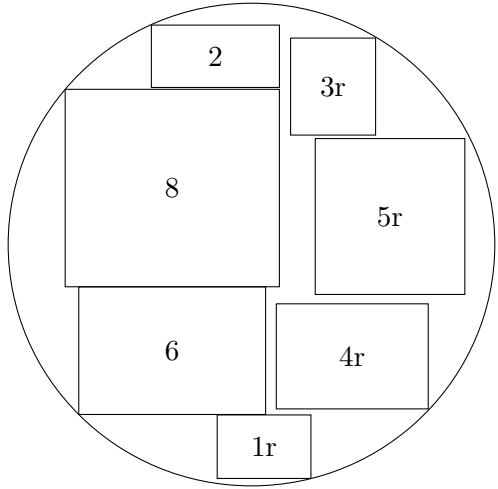
\begin{figure}[!htb]
\centering
\begin{tikzpicture}[scale=0.8];
\draw (0,0) circle [radius=        4.0000000000cm];
\draw(       -1.6459982931 ,         3.6374193994) -- 
(         .4592648648 ,         3.6374193994) -- 
(         .4592648648 ,         2.6039265764) -- 
(       -1.6459982931 ,         2.6039265764) -- 
(       -1.6459982931 ,         3.6374193994) ;
 \node at (        -.5933667142 ,         3.1206729879) {     2     }; 
\draw(       -2.8395332087 ,         -.7018428058) -- 
(         .2322371262 ,         -.7018428058) -- 
(         .2322371262 ,        -2.8166753417) -- 
(       -2.8395332087 ,        -2.8166753417) -- 
(       -2.8395332087 ,         -.7018428058) ;
 \node at (       -1.3036480412 ,        -1.7592590738) {     6     }; 
\draw(       -3.0626262509 ,         2.5729737923) -- 
(         .4589048496 ,         2.5729737923) -- 
(         .4589048496 ,         -.6997534805) -- 
(       -3.0626262509 ,         -.6997534805) -- 
(       -3.0626262509 ,         2.5729737923) ;
 \node at (       -1.3018607006 ,          .9366101559) {     8     }; 
\draw(        -.5635507165 ,        -2.8242854782) -- 
(         .9771191400 ,        -2.8242854782) -- 
(         .9771191400 ,        -3.8769170572) -- 
(        -.5635507165 ,        -3.8769170572) -- 
(        -.5635507165 ,        -2.8242854782) ;
 \node at (         .2067842117 ,        -3.3506012677) {     1r    }; 
\draw(         .6431945584 ,         3.4222022955) -- 
(        2.0403237450 ,         3.4222022955) -- 
(        2.0403237450 ,         1.8145467931) -- 
(         .6431945584 ,         1.8145467931) -- 
(         .6431945584 ,         3.4222022955) ;
 \node at (        1.3417591517 ,         2.6183745443) {     3r    }; 
\draw(         .4079414892 ,         -.9825721177) -- 
(        2.9055491447 ,         -.9825721177) -- 
(        2.9055491447 ,        -2.7241989119) -- 
(         .4079414892 ,        -2.7241989119) -- 
(         .4079414892 ,         -.9825721177) ;
 \node at (        1.6567453169 ,        -1.8533855148) {     4r    }; 
\draw(        1.0480184995 ,         1.7564082241) -- 
(        3.5073486430 ,         1.7564082241) -- 
(        3.5073486430 ,         -.8273238333) -- 
(        1.0480184995 ,         -.8273238333) -- 
(        1.0480184995 ,         1.7564082241) ;
 \node at (        2.2776835713 ,          .4645421954) {     5r    };

\end{tikzpicture}
\captionsetup{font=small}
\caption{Maximise the total area of the rectangles packed, rotation allowed, solution value 37.9687}
\label{fig2a}
\end{figure}

If the rectangles can be rotated through ninety degrees then:
\begin{itemize}

\item If we are wish to maximise the number of rectangles packed Figure~\ref{fig1a} shows the solution as derived by the approach presented in this paper. In that figure we can see that seven of the ten rectangles have been packed, three rectangles are left unpacked.  

\item If we are wish to maximise the total area of the rectangles packed Figure~\ref{fig2a} shows the solution as derived by the approach presented in this paper. In that figure we can see that seven of the ten rectangles have been packed.   
\end{itemize}

In Figure~\ref{fig1a} and Figure~\ref{fig2a} the letter r after the rectangle number indicates that the rectangle has been rotated through ninety degrees. 
Comparing Figure~\ref{fig1} and Figure~\ref{fig1a} we can see that they both involve the packing of seven rectangles.
 Whilst allowing rotation through ninety degrees allows the possibility of a better solution as compared with the no rotation case this is by no means assured.
Comparing Figure~\ref{fig2} and Figure~\ref{fig2a} we can see that in this particular case an improvement in the total area of the rectangles packed has been made by making use of rotation.

The structure of this paper is as follows. 
In 
Section~\ref{Sec:Lit}  we review the literature relating to the packing of rectangles. We discuss application areas where rectangle packing problems arise. 
We also review the literature relating to the particular metaheuristic, formulation space search,  used in this paper. 
In 
Section~\ref{Sec:Formulation} we formulate the problem of packing unequal rectangles/squares into a fixed size circular container as a mixed-integer nonlinear program. We show how we can eliminate a nonlinear maximisation term that arises in one of the constraints in our formulation.  We also show how we can deal with the case where rectangles can be rotated through ninety degrees.
We  
indicate the amendments that can be made to the formulation
for the special case where we are maximising the number of squares packed.
Section~\ref{Sec:FSS} gives details of the formulation space search heuristic that we use to solve the problem. Computational results are  presented in Section~\ref{Sec:Results} for problems involving up to 30 rectangles/squares. In that section we give results both for maximising the number of rectangles/squares packed and for maximising the total area of the rectangles/squares packed.
Finally in Section~\ref{Sec:Conclusions} we present our conclusions.

\section{Literature survey} \label{Sec:Lit}

In this section we first discuss the literature relating to the problem of packing rectangles and its applications.
 We then discuss the literature relating to the particular metaheuristic, formulation space search,  we use to solve the rectangle packing problem considered in this paper.

\subsection{Rectangle packing}

The majority of the work in the literature related to rectangle packing deals with packing rectangles/squares within a larger container that is either a square, or a rectangle, or a rectangular strip with one dimension fixed and the other dimension variable (e.g.~fixed width, but variable length). 

A common feature of such work is that it is assumed that all of the smaller  rectangles have to be packed into the larger container, which leads to an optimisation problem relating  to minimising the dimension of the container. For example for a square container a natural optimisation problem is to minimise the side of the square container (which also minimises its perimeter and area). For a rectangular container one can examine minimising either its perimeter or its area. For a rectangular strip one can minimise the variable dimension.
With respect to the packing of rectangles within a circular container then
the natural optimisation problem is to minimise the radius of the container.

In our literature survey below we focus principally on papers that take a packing approach. The reader may be aware that a closely related problem to packing is cutting e.g.~cutting rectangles from a larger stock rectangle. There has been a substantial amount of work presented in the literature dealing with cutting. However much of that work involves additional restrictions with regard to the cuts that are made. One such restriction might be that the cuts are guillotine cuts, a guillotine cut on a rectangle being a cut from one edge of the rectangle to the opposite edge which is parallel to the two remaining edges. Another such restriction might be to limit the cutting to a number of stages, where at each stage guillotine cuts are made, but in a direction opposite to that adopted in the previous stage. So for example in the first stage guillotine cuts are made parallel to the $y$-axis, then in the second stage guillotine cuts are made parallel to the $x$-axis, etc. Since the primary focus of the work presented in this paper is packing rectangles within a 
\emph{circular container} we, for space reasons, exclude detailed consideration of work focused on cutting rectangles from rectangular containers from the literature survey presented below.

Unless otherwise stated all of the work considered below deals with orthogonal packing, so rectangles/squares are packed without rotation.

Li and Cheng~\cite{Li89} show that the problem of determining
whether a set of squares can be packed into a larger rectangle
is strongly NP-complete. In addition they show that the problem
of determining
whether a set of
rectangles can be packed  into a square is NP-complete.
Leung et al.~\cite{Leung1990} show that the problem of determining
whether a set of squares can be  packed 
into a square is strongly NP-complete. 

Picouleau~\cite{Picouleau1996} considered the worst-case analysis of three fast heuristics for packing squares into a square container so as to minimise the size of the square.
Murata et al.~\cite{Murata1996} present a simulated annealing algorithm for the problem for  packing rectangles into a  rectangular container so as to minimise the size (area) of the container.
Liu and Teng~\cite{Liu1999} present a genetic algorithm for the problem of packing a set of rectangles into a strip of fixed width using minimum height.

Wu et al.~\cite{Wu2002} present a heuristic attempting to pack every member of a  set of rectangles inside a fixed size rectangular container.
Caprara et al.~\cite{Caprara2006} discuss absolute worst-case performance ratios for lower bounds on packing rectangles/squares into a  square container so as to minimise the size of the square container. They consider the case where the rectangles have fixed orientation and the case where they can be rotated through ninety degrees.
Huang et al.~\cite{Huang2007}  present a heuristic approach to packing rectangles within a fixed size rectangular container so as to maximise the total area of the rectangles packed where the rectangles can be rotated through ninety degrees.

Birgin et al.~\cite{Birgin2010a} consider packing the maximal number of identically sized rectangles inside a rectangular container. Their approach is based upon recursive partitioning and allows the rectangles to be rotated through ninety degrees.
Korf et al.~\cite{Korf2010} consider the problem of packing a set of rectangles (with and without ninety degree rotation allowed) in a rectangular container of minimal area. They adopt a constraint satisfaction approach to the problem.
Maag et al.~\cite{Maag2010} consider the problem of packing a set of rectangles  in a rectangular container of minimal area. Their approach is based on relaxing the constraint on rectangle overlap.

Huang and Korf~\cite{Huang2012}  consider the same problem 
as~\cite{Korf2010} but adopt an approach based on  first deciding 
$x$-coordinate values for each rectangle. 
Bortfeldt~\cite{Bortfeldt2013} presents a number of heuristic approaches (based on solution methods for two-dimension knapsack and two-dimension strip packing) for  packing rectangles into a  rectangular container so as to minimise the size (area) of the container.

Martello and Monaci~\cite{Martello2015} consider 
the problem of packing rectangles/squares into a  square container so as to minimise the size of the container. They present a linear integer programming formulation and an exact approach based on a 
two-dimensional packing algorithm as well as a  metaheuristic. 
They deal with the case where the rectangles have fixed orientation and also the case where they can be rotated through ninety degrees.
Delorme et al.~\cite{Delorme2017} present a Benders' decomposition approach to the problem of packing a set of rectangles (with ninety degree rotation allowed) into a strip of fixed width using minimum height. Their approach (as they discuss) can be easily applied to the problem of packing rectangles/squares into a  square container of minimal size. 

It is important to note here that a number of the approaches given in the literature for the problem of packing rectangles within a rectangular container utilise the fact that rectangle position coordinates can be taken from a finite discrete set (e.g.~by packing rectangles so that they are positioned at their lowest bottom-left position). For example 
see~\cite{Delorme2017,Martello2015}.
However in this paper we consider a circular container, and 
\emph{the lack of rectangular sides to the container render such discretisation approaches invalid for the problem we consider}.

As far as we are aware the problem considered in this paper of packing unequal rectangles/squares into a circular container has only been considered by just a few papers in the literature previously.
Li et al.~\cite{Li2014} consider the problem of packing orthogonal unequal rectangles in a circular container with an additional constraint related to mass balance. Their objective function is to minimise the radius of the container. A heuristic algorithm is presented.

Hinostroza et al.~\cite{Hinostroza2013} consider the problem of cutting rectangular boards from a log, regarded as a circular container. They present a nonlinear formulation of the problem (based on~\cite{Birgin2006}), and two heuristics, one based on ordering the rectangles and the other on simulated annealing. Note here that, in our judgement, their formulation is flawed.

Work has been presented in the literature relating to packing rectangles/squares into arbitrary convex regions, and such work can be applied to a circular container. We discuss this work below.

Birgin at al~\cite{Birgin2006a} introduce the concept of sentinels sets, which are finite subsets of the items to be packed such that, when two
items are superposed, at least one sentinel of one item is in the interior of the other item. Using these sentinel sets they consider packing identical rectangles within both convex regions and a rectangular container, with and without rectangle rotation (both ninety degree rotation and arbitrary rotation).
Birgin et al.~\cite{Birgin2006} consider packing rectangles (with and without ninety degree rotation). Their objective is to feasibly pack all rectangles. Iteratively increasing the number of rectangles enables one to maximise the number of (identical) rectangles placed.  Their approach is based on nonlinear optimisation.

Birgin and Lobato~\cite{Birgin2010} consider 
packing identical rectangles within an arbitrary convex region where a common rotation of $\theta$ degrees (not restricted to $\theta=90$) of all the rectangles is allowed. In addition a rectangle can be rotated through ninety degrees before a rotation of $\theta$ is applied. 
 Their solution method is a combination
of branch and bound and active-set strategies for bound-constrained minimization of smooth
functions. 
Cassioli and Locatelli~\cite{Cassioli2011} present a heuristic approach based on iterated local search for the problem of packing the maximum number of rectangles of the same size  within a convex region (where rectangle rotation through ninety degrees is allowed).

Andrade and Birgin~\cite{Andrade2013} present symmetry breaking constraints for two problems relating to packing identical rectangles (with or without ninety degree rotation) in a polyhedron. They consider packing as many identical rectangles as possible within a given polyhedron as well as finding the smallest polyhedron of a specified type that accommodates a fixed number of identical rectangles.

More generally Birgin~\cite{Birgin2016} considers the application of nonlinear programming in packing problems. They note that nonlinear
programming formulations and methods have been successfully applied to a wide range of packing problems. In particular we in this paper, as in the formulation presented below, use a nonlinear model.

\subsection{Applications}

The problem of packing  rectangular objects into a larger container (equivalently cutting rectangular objects from a larger container)
appears in a number of practical situations. As noted in Dowsland and Dowsland~\cite{Dowsland1992} the earliest applications were in  glass and metal industries where smaller rectangular objects had to be cut from larger (typically rectangular) stock pieces. A further application they discuss occurs in pallet loading where rectangular boxes have to be packed onto a wooden pallet for transport. 
Sweeney and  Paternoster~\cite{Sweeney1992} present an 
application-orientated research bibliography that lists some of the early work related to packing. 

Lodi et al.~\cite{lodi2002} present a literature survey relating to two-dimensional packing and solution approaches. They mention a number of practical applications relating to rectangle cutting/packing. These include the arrangement of articles and advertisements on newspaper pages and   in the wood and glass industry cutting rectangular
items from larger sheets of material. They also mention the placement of goods on shelves in warehouses.
Wascher et al.~\cite{ wascher2007} also mention some practical applications (such as pallet loading) in their work presenting  a typology of cutting and packing problems.

In relation to the specific problem considered in this paper of packing unequal rectangles/squares into a fixed size circular container we are aware of a number of practical applications. 

For example in the forestry/lumber industry consider the cutting of rectangular wooden boards from timber logs made from trees that have been felled.
Here, by approximating the shape of the timber log by a circle of known radius, we have the problem considered in this paper,
namely which of the rectangles (of known sizes) that we desire to cut should be cut from the circular log~\cite{Hinostroza2013}.

A further practical example relates to the problem considered 
in~\cite{Li2014} which was concerned with packing orthogonal unequal rectangles in a circular container with an additional constraint related to mass balance. Here the container was a satellite and the rectangular objects related to items comprising the satellite payload. The mass balance constraint considered in~\cite{Li2014} was a single nonlinear constraint that involved the (mass weighted) centres of each rectangle. Since, as will become apparent below, our formulation space search approach for packing rectangles into a fixed size circular container is based on a mixed-integer nonlinear program  it is trivial to introduce into our approach a single additional nonlinear constraint (such as a mass balance constraint).

\subsection{Formulation space search}
When solving nonlinear non-convex problems with the aid of a solver, 
Mladenovi\'{c} et al.~\cite{Mladenovic2005} observed that different formulations of 
the same problem may have different characteristics.  
Hence a natural way to proceed is by
swapping between formulations. Under this framework Mladenovi\'{c} 
et al.~\cite{Mladenovic2005} use formulation space search (henceforth FSS) for the circle packing problem considering two formulations of the problem: one in a 
Cartesian coordinate system, the other in a Polar coordinate system. Their algorithm solves the problem with one formulation 
at a time and when the solution is the same for all formulations  the algorithm terminates. They consider
packing identical circles 
into the unit circle and the unit square. 

In Mladenovi\'{c} et al.~\cite{Mladenovic2007} they improve 
on~\cite{Mladenovic2005} by considering a mixed 
formulation of the problem. They set a subset of the circles in the Cartesian system whilst the rest of 
the circles were in the Polar system. 
López and Beasley~\cite{BL2011} use FSS for the problem of packing 
equally sized circles inside a variety of containers. They present computational results which 
show that their approach improves upon previous results based on FSS presented in the literature.
For some of the containers considered they improve on the best result previously known.
López and Beasley~\cite{BL2013} use FSS  to solve the packing problem with non-identical circles in different shaped containers. They present computational results which were compared with benchmark problems and also proposed some new instances.
López and Beasley~\cite{BL2016} use FSS  to solve the problem of packing  non-identical circles in a fixed size container. 

Essentially FSS exploits the fact that:
\begin{compactitem}
\item because of the nature of the solution process in nonlinear optimisation we often fail to obtain a globally optimum solution from a single formulation; and so
\item perturbing/changing the formulation and then resolving the nonlinear program may lead to an improved solution. 
\end{compactitem}
Given the above it is a simple matter to construct iterative schemes that move between formulations in a systematic manner.

FSS has been applied to a few problems  additional to circle packing (e.g.~timetabling~\cite{ Kochetov2008}). In~\cite{BL2014}  FSS was used   to solve some benchmark mixed-integer nonlinear programming problems. In a more general sense an adaptation to  FSS was presented in~\cite{Brimberg2014} for solving continuous location problems. More discussion as to FSS can be found in Hansen et al.~\cite{Hansen2010}. A related approach is variable space search, which has been applied to graph colouring (Hertz et al.~\cite{Hertz2008, Hertz2009}). Other related approaches are variable formulation search which has been applied to the cutwidth minimisation 
problem~\cite{Pardo2013,Duarte2016}
and variable objective search which has been applied to  the maximum independent set problem~\cite{butenko13}.

As noted in Pardo et al.~\cite{Pardo2013}  
variable space search, 
variable formulation search and 
variable objective search 
contain similar ideas as originally expounded using FSS. At a slightly more general level FSS can be regarded as a variant of variable neighbourhood search, for example see~\cite{Amirgaliyeva2017, Hansen2017}.

\section{Formulation} \label{Sec:Formulation}

In this section we first present our basic formulation for the problem of packing unequal rectangles in a fixed size circular container as a mixed-integer nonlinear program (MINLP).  We then show how we can eliminate a nonlinear maximisation term that arises in one of the constraints in our formulation. We indicate how we can deal with the case where rectangles can be rotated through ninety degrees.
For the special case where we are maximising the number of squares packed we present the amendments that can be made to the formulation.

\subsection{Basic formulation} 

The problem we consider is to find the maximal weighted packing of 
$n$ unequal rectangles in a fixed size circular container. Here we have the option, for each unequal rectangle, of choosing to pack it or not. We can formulate this problem as follows.

Let the fixed size circular container be 
 of radius $R$ and, without loss of generality, let it 
be centred at the origin of the Euclidean
plane. We have $n$ rectangles from which to construct a packing, where rectangle $i$ has a horizontal side of length $L_i$ and a vertical side of width $W_i$,  and value (if packed) $V_i$. In 
our basic  formulation we do not allow any rotation when packing rectangles so that rectangles are packed with their horizontal (length) edges parallel to the $x$-axis, their vertical (width) edges parallel to the $y$-axis.
Clearly if we are dealing with packing squares then $L_i=W_i$.
Here we label the rectangles so that they are ordered in increasing size (area) order 
(i.e.~$L_iW_i \leq L_{i+1}W_{i+1}~i=1,\ldots,n-1$).

Using a value $V_i$ here for each rectangle $i$ enables us to consider a number of different problems within the same formulation. For example if we take $V_i = 1~i=1,\ldots,n$ then we have the problem of maximising the number of rectangles packed. If we take $V_i =  L_iW_i~i=1,\ldots,n$ then we have the problem of maximising the total area of the rectangles packed. Alternatively the $V_i~i=1,\ldots,n$ can be assigned arbitrary values.

Then the variables are:
\begin{compactitem} 
\item $\alpha_i = 1$ if rectangle $i$ is packed, 0 otherwise;~$i=1, \ldots, n$
\item $(x_i,y_i)$ the position of the centre of rectangle $i$;~$i=1, \ldots, n$
\end{compactitem} 
\noindent With regard to the positioning (so $(x_i,y_i)$) of any unpacked rectangle $i$ (for which $\alpha_i=0$) our formulation forces all unpacked rectangles to be positioned at the origin. Let $Q$ be the set of all rectangle pairs $[(i,j) \ | \ i=1,...,n; \ j=1,...,n; \ j > i ]$. The formulation is:
\begin{align}
\max & \hspace{.3cm} \sum_{i=1}^{n} \alpha_iV_i & & & \label{e1} \\ \notag
\text{subject to}& & & & \\
& -\alpha_i(\sqrt{(R^2-W_i^2/4)}-L_i/2) \leq x_i \leq \alpha_i(\sqrt{(R^2-W_i^2/4)}-L_i/2) & & i=1,\ldots,n \label{e4} \\
 & -\alpha_i(\sqrt{(R^2-L_i^2/4)}-W_i/2) \leq y_i \leq \alpha_i(\sqrt{(R^2-L_i^2/4)}-W_i/2) & & i=1,\ldots,n \label{e5} \\
 & (x_i+L_i/2)^2 + (y_i+W_i/2)^2 \leq  \alpha_iR^2 
+(1-\alpha_i)(L_i^2/4+W_i^2/4)  & & i=1,\ldots,n \label{e2a} \\
 & (x_i+L_i/2)^2 + (y_i-W_i/2)^2 \leq  \alpha_iR^2 
+(1-\alpha_i)(L_i^2/4+W_i^2/4)  & & i=1,\ldots,n \label{e2b} \\
 & (x_i-L_i/2)^2 + (y_i+W_i/2)^2 \leq  \alpha_iR^2 
+(1-\alpha_i)(L_i^2/4+W_i^2/4) & & i=1,\ldots,n \label{e2c} \\
 & (x_i-L_i/2)^2 + (y_i-W_i/2)^2 \leq  \alpha_iR^2 
+(1-\alpha_i)(L_i^2/4+W_i^2/4) & & i=1,\ldots,n \label{e2d} \\
& \alpha_i\alpha_j[\max\{|x_i-x_j| - (L_i+L_j)/2, |y_i-y_j| - (W_i+W_j)/2\}]
 \geq 0 & & \forall (i,j) \in Q
\label{e3} \\
& \alpha_i \in \{0,1\}
& & i=1,\ldots,n \label{e6} 
\end{align}

The objective function, Equation~(\ref{e1}), maximises the value of the rectangles packed. 
Equation~(\ref{e4}) ensures that if a rectangle is packed 
(i.e.~$\alpha_i=1$) its 
$x$-coordinate lies in $[-(\sqrt{(R^2-W_i^2/4)}-L_i/2),+(\sqrt{(R^2-W_i^2/4)}-L_i/2)]$. These limits can be easily deduced from geometric considerations, e.g.~consider the centre $x$-coordinate value associated with a rectangle placed with its centre on the 
$x$-axis and with two of its corners just touching the circular container.
The key feature of Equation~(\ref{e4}) is that if the rectangle is not packed (i.e.~$\alpha_i=0$) then the $x$-coordinate is forced to be zero. Equation~(\ref{e5}) is the equivalent constraint to 
Equation~(\ref{e4})
for the $y$-coordinate.

Equations~(\ref{e2a})-(\ref{e2d}) ensure that if a rectangle is packed (so for rectangle $i$ with $\alpha_i=1$) its centre is appropriately positioned such that the entire rectangle lies inside the circular container. To achieve this we need to ensure that all four corners of the rectangle lie inside the circular container. These four corners are $(x_i \pm L_i/2, y_i \pm W_i/2)$ and 
Equations~(\ref{e2a})-(\ref{e2d}) ensure that the 
(squared)
distance from the origin to each these corners is no more than the 
(squared)
radius of the container. Note that if the rectangle is packed (so $\alpha_i=1$) the left-hand side of Equations~(\ref{e2a})-(\ref{e2d}) is $R^2$.

If the rectangle is not packed (so $\alpha_i=0$) then from Equations~(\ref{e4}),(\ref{e5}) the rectangle is positioned at the origin (so has $x_i=y_i=0$). In that case the left-hand side of Equations~(\ref{e2a})-(\ref{e2d}) becomes $L_i^2/4 + W_i^2/4$, 
as does the right-hand side, and so the constraints are automatically satisfied.

Equation~(\ref{e3}) guarantees that any two rectangles $i$ and $j$ which are both packed (so $\alpha_i=\alpha_j=1$) do not overlap each other. This constraint is derived from that given previous 
in~\cite{Chr74}. 
It states that two rectangles of size $[L_i,W_i]$ and $[L_j,W_j]$ do not overlap provided that 
the difference between their centre $x$-coordinates is at least $(L_i+L_j)/2$ or that 
the difference between their centre $y$-coordinates is at least $(W_i+W_j)/2$ (or both). If one or other of the rectangles is not packed the left-hand side of Equation~(\ref{e3}) becomes zero due to the product term ($\alpha_i\alpha_j$) which means that the constraint is automatically satisfied. Equation~(\ref{e6}) is the integrality constraint. 

As discussed above our formulation positions any unpacked rectangle at the origin. For unpacked rectangle 
$i$ the inclusion of an appropriate $\alpha_i$ term on the left-hand side of Equation~(\ref{e3}) ensures that this unpacked rectangle, although positioned at the origin, does not actively participate in the overlap constraint which must apply between all packed rectangles. 

Our formulation (Equations~(\ref{e1})-(\ref{e6})) is a mixed-integer nonlinear program (MINLP). Computationally MINLPs are recognised to be very demanding, involving as they do both an element of combinatorial choice and solution of an underlying  continuous nonlinear program. For the problem considered in this paper the combinatorial choice relates to the choice of the set of rectangles to be packed, and the underlying  continuous nonlinear program relates to deciding where to feasibly position within the circular container the rectangles that are packed.

\subsection{Elimination of the maximisation term} 

The overlap constraint (Equation~(\ref{e3})) contains the  expression
$\max\{|x_i-x_j| - (L_i+L_j)/2, |y_i-y_j| - (W_i+W_j)/2\}$.  
 For the particular problem considered in this paper this maximisation term  can be eliminated, 
albeit by enlarging the size of the MINLP to be solved.

Introduce additional continuous variables 
$\beta_{ij},~\forall (i,j) \in Q,$ 
defined by:
\begin{align}
& 0 \leq \beta_{ij} \leq 1  & & \forall (i,j) \in Q
\label{eb1} 
\end{align}

Then we  can replace  Equation~(\ref{e3})   by:
\begin{align}
& \alpha_i\alpha_j \big[
\beta_{ij}[|x_i-x_j| - (L_i+L_j)/2] + (1-\beta_{ij})[ |y_i-y_j| - (W_i+W_j)/2] \big]
 \geq 0 & & \forall (i,j) \in Q
\label{eb2} 
\end{align}

The logic here is that
the $\alpha_i\alpha_j$ term
ensures that the Equation~(\ref{eb2}) is always satisfied when either $\alpha_i=0$ or $\alpha_j=0$ (as indeed it does in Equation~(\ref{e3})).
It only remains to check therefore the validity of replacing
Equation~(\ref{e3}) with Equation~(\ref{eb2})
in the case 
$\alpha_i=\alpha_j=1$.

When $\alpha_i=\alpha_j=1$ Equation~(\ref{eb2}) becomes
$\beta_{ij}[|x_i-x_j| - (L_i+L_j)/2] + (1-\beta_{ij})[ |y_i-y_j| - (W_i+W_j)/2] \geq 0 $. 
Now the weighted sum on the left-hand side of this constraint
 \emph{can only be non-negative provided that at least one of the two terms in it is itself non-negative}. In other words Equation~(\ref{eb2}) will ensure that one (or both) of $[|x_i-x_j| - (L_i+L_j)/2]$ and $[ |y_i-y_j| - (W_i+W_j)/2]$ will be 
non-negative. Since one or both of these terms are non-negative it is therefore true that the maximisation term in
Equation~(\ref{e3}), 
$\max\{|x_i-x_j| - (L_i+L_j)/2, |y_i-y_j| - (W_i+W_j)/2\}$,
 must  also be non-negative. This in turn implies that Equation~(\ref{e3}) is satisfied. Therefore it is valid to replace Equation~(\ref{e3}) by Equation~(\ref{eb2}).

Note here that it is also valid 
 to replace Equation~(\ref{e3}) by Equation~(\ref{eb2})
if we define $\beta_{ij}$ as binary (zero-one) variables. However we might well expect there to be computational benefit in defining these variables as continuous, rather than binary, variables.

\subsection{Rotation}
As is common in the literature 
(e.g.~\cite{Birgin2010a,  Caprara2006, Hinostroza2013, Huang2012, Korf2010,  lodi2002, Maag2010,  Martello2015}) in the basic formulation presented above we did not allow any rotation when packing,  so that the items to be packed (rectangles/squares) were packed with their horizontal (length) edges parallel to the $x$-axis, their vertical (width) edges parallel to the $y$-axis.
If rotation of any item is allowed (which might be dependent on the practical problem being modelled) then the situation becomes more complex, although obviously rotation might enable a better solution to be found.

In  the literature rotation through ninety degrees is the most common situation modelled (e.g.~\cite{ Birgin2010a,  Caprara2006,   Delorme2017,   Huang2007,  Huang2012,  Korf2010,  Li2014,  Martello2015,  Murata1996, Wu2002}). 
Clearly rotation through ninety degrees is  irrelevant when we are packing squares (as they are the same under ninety degree rotation) and  only relevant when we are dealing with unequally sized rectangles. Our formulation can be extended to deal with rotation through ninety degrees as discussed below. Rotation through an
arbitrary 
 angle cannot be dealt with by our approach.

If the rectangles can be rotated through ninety degrees then this is easily incorporated into our formulation. Suppose that rectangle $i$ can be rotated through ninety degrees. Then create a new rectangle ($j$ say) that represents rectangle $i$ if it is rotated, so that we have $L_j=W_i$,~$W_j=L_i$,~$V_j=V_i$. Add to the formulation:
\begin{align} 
& \alpha_i + \alpha_j \leq 1 \label{rot1}
\end{align}
Equation~(\ref{rot1}) ensures that we cannot use both the original rectangle $i$ and its rotated equivalent $j$. Dealing with rectangle rotation therefore requires creating a new rectangle for each original rectangle that can be rotated and adding a single constraint to the formulation.

In terms of the effect on the formulation then if all $n$ rectangles can be rotated this only directly adds $n$ constraints (Equation~(\ref{rot1})) to the formulation. However the creation of an additional $n$ rotated rectangles doubles the number of rectangles to be considered for packing. This means that the number of linear
constraints associated with Equations~(\ref{e4}),(\ref{e5}) doubles, as does the number of nonlinear constraints associated with Equations~(\ref{e2a})-(\ref{e2d}). The more significant effect is that  the number of nonlinear
constraints associated with Equation~(\ref{e3}) increases from 
$n(n-1)/2$ to $2n(2n-1)/2$ (so approximately increases by a factor of 4). This increase in the number of nonlinear
constraints associated with Equation~(\ref{e3}) also carries through to increase the number of $\beta_{ij}$ variables (Equation~(\ref{eb1}))  that need to be considered by an (approximate) factor of 4. For this reason we would expect that, computationally, dealing with a problem with $n$ rectangles with fixed orientation becomes much more challenging if all $n$ rectangles can be rotated.

\subsection{Maximising the number of squares packed} \label{Sec:opt}

From our previous work~\cite{BL2016} we know that when we are considering a packing problem where all the items to be packed can be ordered such that item $i$ fits inside item $j$ for all $j>i$ then, in the case where we are maximising the number of items packed, 
  the optimal solution consists of the first $K$ items, 
for some $K$.

Clearly items can be ordered to fit inside each other if we are considering packing squares, i.e.~order the squares in increasing size (length) order, but such an ordering is unlikely to be possible if we are packing rectangles. Hence we shall just consider square packing here.

In the case of square packing therefore, when we are maximising the number of squares packed, 
 we can impose the additional constraints:
\begin{align} 
& \alpha_{i-1} \geq \alpha_i & i=2,\ldots,n \label{e9} \\ 
& \alpha_k = 0 & \text{if}~\sum_{i=1}^{k} L_i^2 > \pi R^2 & & k=1,\ldots,n \label{e10}
\end{align}
Equation~(\ref{e9}) ensures that if $\alpha_i$ is one (so square $i$ is packed) then $\alpha_{i-1}$ must also be one (so square $i-1$ is packed). If square $i$ is not packed ($\alpha_i=0$) then the right-hand side of this constraint is zero, so the constraint is always satisfied whatever the value for $\alpha_{i-1}$. 
Collectively the $(n-1)$ inequalities represented in 
Equation~(\ref{e9}) ensure that  the optimal solution consists of the first $K$ squares, for some $K$.

In Equation~(\ref{e10}) we have that
if we were to pack square $k$ then we would have to pack all squares up to and including square $k$. If this packing exceeds the area of the container then clearly square $k$ cannot be packed.

Aside from these additional constraints we can  amend the overlap constraint, Equation~(\ref{eb2}). Note that Equation~(\ref{eb2})  includes a $\alpha_i\alpha_j$ term and applies for $(i,j) \in Q$, where $Q$ is defined to have $j>i$.
Now if $\alpha_j=1$ we automatically know that $\alpha_i=1$ (since $j>i$) and hence that $\alpha_i\alpha_j=1$. If $\alpha_j=0$ then it is irrelevant what value $\alpha_i$ takes since we must have $\alpha_i\alpha_j=0$. In other words the $\alpha_i\alpha_j$ term in Equation~(\ref{eb2}) can be replaced by $\alpha_j$ so that the overlap constraint, Equation~(\ref{eb2}), becomes:
\begin{align}
& \alpha_j \big[ \beta_{ij}[|x_i-x_j| - (L_i+L_j)/2] + (1-\beta_{ij})[ |y_i-y_j| - (W_i+W_j)/2] \big] 
 \geq 0 & & \forall (i,j) \in Q
\label{eb2a} 
\end{align}
Note here that we have used $W_i$ and $W_j$ in Equation~(\ref{eb2a}) for clarity of comparison with Equation~(\ref{eb2}). Obviously since we are just considering square packing here we have 
$W_i=L_i~i=1,\ldots,n$.

\section{FSS algorithm} \label{Sec:FSS}

In this section we  present our FSS algorithm for the problem. For simplicity we  present our approach using the basic formulation of the problem, Equations~(\ref{e1})-(\ref{e6}), before the amendments as discussed above (i.e.~elimination of the maximisation term and adaptions for packing squares).
We discuss at the end of this section how we incorporate a number of other constraints, presented in this section, into our approach.

\subsection{Algorithm}
Consider the formulation, Equations~(\ref{e1})-(\ref{e6}), given above. Letting $\delta$ be a small positive constant  replace the integrality requirement, Equation~(\ref{e6}), by:
\begin{align}
& \sum_{i=1}^{n} \alpha_i(1-\alpha_i) \leq \delta \label{e7} \\ 
& 0 \leq \alpha_i \leq 1 & i=1,\ldots,n \label{e8} 
\end{align}

If $\delta$ was zero these equations would force $[\alpha_i,~i=1,\ldots,n]$ to assume zero-one values. 
However given the capabilities of nonlinear optimisation software 
simply replacing an explicit integrality condition by 
Equations~(\ref{e7}),(\ref{e8})
 would not 
be computationally successful, since we would be hoping to generate a (globally optimal) solution to a 
continuous nonlinear optimisation problem with a very tight inequality constraint. Note that if $\delta$ is zero then 
Equation~(\ref{e7}) is effectively an equality constraint as the left-hand side is non-negative. 

Accordingly we  adopt a heuristic approach and have $\delta > 0$. Hence our original MINLP,
Equations~(\ref{e1})-(\ref{e6}),
 has now become a continuous nonlinear optimisation problem, since we have relaxed the integrality requirement using Equations~(\ref{e7}),(\ref{e8}). This nonlinear optimisation problem is optimise Equation~(\ref{e1}) subject to 
Equations~(\ref{e4})-(\ref{e3}),(\ref{e7}),(\ref{e8}). We refer to this problem as the \emph{continuous FSS relaxation} of the problem.

If we solve this nonlinear problem the $[\alpha_i,~i=1,\ldots,n]$ can deviate (albeit only slightly, if $\delta$ is small) from their ideal zero-one values, but we can round them to their nearest integer value to recover an integer set of values. Given an integer set of values for $[\alpha_i,~i=1,\ldots,n]$ then the original formulation (Equations~(\ref{e1})-(\ref{e6})) becomes a nonlinear feasibility problem. This 
nonlinear feasibility problem is to
 find positions $(x_i,y_i)$ for each rectangle $i$ that we have chosen to pack (so with $\alpha_i=1$ in the rounded solution). Note here that this is a feasibility problem as the objective function, Equation~(\ref{e1}), is purely a function of the zero-one variables (and these have been fixed by rounding).

With just a single value for $\delta$ we have just a single nonlinear problem: optimise Equation~(\ref{e1}) subject to
Equations~(\ref{e4})-(\ref{e3}),(\ref{e7}),(\ref{e8}).
 However changing $\delta$ is a systematic fashion creates a series of different problems that can be given to an appropriate nonlinear solver in an attempt to generate new and improved solutions to our original MINLP. The idea here is that altering $\delta$ perturbs the nonlinear formulation and hence, given the nature of any nonlinear solution software, might lead to a different solution.

The pseudocode for our FSS algorithm for the rectangle packing problem considered in this paper is presented in Algorithm~\ref{fssalg}. In this pseudocode let $P$ denote the original MINLP (here optimise Equation~(\ref{e1}) subject to
Equations~(\ref{e4})-(\ref{e6})) and $P^*$ denote the continuous FSS relaxation (here optimise Equation~(\ref{e1}) subject to
Equations~(\ref{e4})-(\ref{e3}),(\ref{e7}),(\ref{e8})).

We first initialise values, here $Z_{best}$ is the best feasible solution found and $t$ is an iteration counter. We then solve the original MINLP $P$. 

In this pseudocode all attempts to solve a nonlinear problem (e.g.~$P$ or $P^*$) must be subject to a time limit, since otherwise the computation time consumed could become extremely high. For this reason we always terminate the solution process after a predefined time limit, returning the best feasible solution found (if one has been found). In the computational results reported later below this time limit was set to $10n$ seconds. 

SCIP is capable of solving our MINLP formulation $P$ 
 to proven global optimality because SCIP restricts the type of nonlinear expression allowed~\cite{bussieck14, vigerske17, vigerske16}. However with regard to our computational results for all the problem instances considered below this never occurred within the time limit imposed.

Note here that even if solving $P$ or $P^*$ returns a feasible solution we have no guarantee that this is an optimal solution, since a better solution might have been found had we increased the time limit.

The iterative process in the pseudocode is to update the iteration counter and solve the continuous FSS relaxation $P^*$. 
If  a feasible solution for $P^*$ has been found then we round that solution and 
solve the resulting feasibility problem (again subject to the predefined time limit). The best feasible solution found (if any)
is updated and provided we have not reached the termination condition we reduce $\delta$ by a factor $\gamma$ and repeat.

We terminate when $\delta$ is small ($\leq 10^{-5}$) or we have performed a number of consecutive iterations (three iterations) without improving the value of the best solution found. 
We reduce $\delta$ by a factor $\gamma=0.5$ at each iteration and replicate (repeat) our heuristic a number of times (five replications were performed in the computational results reported below). 
The values for these factors were set based on our previous computational experience with FSS.

\begin{algorithm}[!htb]
\phantom{a}
\caption[FSS pseudocode ]{Formulation space search pseudocode}
\begin{algorithmic}
\STATE \emph{Initialisation:} \ $\delta \leftarrow 0.05$ \ $Z_{best} \leftarrow -\infty$ \ $t \leftarrow 0$ \\
\STATE Solve $P$ and update $Z_{best}$ if a feasible solution for $P$ has been found \\
\STATE \emph{Iterative process:} \\
\WHILE{not termination condition} 
\STATE Update the iteration counter $t \leftarrow t+1$ \\
\STATE Solve $P^*$ \\
\IF{a feasible solution for $P^*$ has been found} 
\STATE Round the $[\alpha_i,~i=1,\ldots,n]$ values in the $P^*$ solution and solve the resulting feasibility problem \\
\STATE Update $Z_{best}$ using the solution to the feasibility problem if a feasible solution for that problem has been found  \\
\ENDIF
\STATE If \ $\delta \leq 10^{-5}$ or $Z_{best}$ has not improved in the last three iterations stop \\
\STATE Update $\delta \leftarrow \gamma\delta$ \\
\ENDWHILE
\end{algorithmic}
\label{fssalg}
\end{algorithm}

\subsection{Constraints}

There are a number of general constraints that apply whatever the objective adopted. Recall that $Z_{best}$ is the value of the best feasible solution encountered during our FSS heuristic. Let the set of rectangles that are packed in this best feasible solution be denoted by $F$. Then the general  constraints that apply are:
\begin{align}
& \alpha_i + \alpha_j \leq 1 & \forall (i,j) \in Q
\hspace{.2cm} \textrm{min}(L_i+L_j,W_i+W_j)>2R \label{e12} \\
& \sum_{i=1}^{n} \alpha_i  L_iW_i \leq \pi R^2 \label{e13} \\
 & \sum_{i=1}^{n} \alpha_iV_i \geq Z_{best} \label{e14} \\
& \sum_{i \not\in F} \alpha_i + \sum_{i \in F} (1-\alpha_i) \geq 1 \label{e15}
\end{align}

Equation~(\ref{e12}) says that if the minimum of the sum of the sides of any two rectangles is greater than the container diameter then we cannot pack both rectangles. Equation~(\ref{e13}) ensures that the total area of the rectangles packed cannot exceed the area of the container. 
Equation~(\ref{e14}) ensures that the value of any solution found is at least that of the best feasible solution known. Equation~(\ref{e15}) is a feasible solution exclusion constraint and ensures that whatever solution is found must differ from the best known solution (of value $Z_{best}$ with packed rectangles $F$) by at least one rectangle. The effect of  Equations~(\ref{e14}) and (\ref{e15}) is to seek an improved feasible solution. 

Note here that although these constraints may be redundant in the original MINLP they may not be redundant in any relaxation of the problem, in particular here the continuous FSS relaxation when we drop the requirement that the $[\alpha_i,~i=1,\ldots,n]$ are zero-one.

\subsection{Summary}
We have presented a considerable number of constraints above and so here (for clarity) we specify the constraints that are involved with $P$ (the original MINLP) and $P^*$ (the continuous FSS relaxation) that are used in the statement of our FSS heuristic given above 
(see Algorithm~\ref{fssalg}). 
\begin{itemize}
\item $P$ is optimise Equation~(\ref{e1}) subject to 
Equations~(\ref{e4})-(\ref{e2d}),(\ref{e6})-(\ref{eb2}),(\ref{e12})-(\ref{e15})
\item $P^*$ is optimise Equation~(\ref{e1}) subject to 
Equations~(\ref{e4})-(\ref{e2d}),(\ref{eb1}),(\ref{eb2}),(\ref{e7})-(\ref{e15})
\end{itemize}

When considering just the packing of squares, so as to maximise the number of squares packed, we add Equations~(\ref{e9}),(\ref{e10}) to   $P$ and  $P^*$ and replace Equation~(\ref{eb2}) by Equation~(\ref{eb2a}). When rectangles can be rotated through ninety degrees we amend the problem in the manner discussed above when we considered Equation~(\ref{rot1}).

\section{Results} \label{Sec:Results}

The computational results presented below 
(Windows  2.50GHz pc, Intel i5-2400S processor, 6Gb memory)
are for our formulation space search heuristic as coded in FORTRAN. We used SCIP (Solving Constraint Integer Programs, version 
4.0.1)~\cite{Achterberg2009,maher17,scip} as the 
mixed-integer nonlinear solver. 
For a technical explanation as to how SCIP solves MINLPs see 
Vigerske and Gleixner~\cite{vigerske16}.
To  input our formulation into SCIP we made use of the modelling language ZIMPL~(Zuse Institute Mathematical Programming Language), and to solve continuous nonlinear problems we used Ipopt~(Interior Point OPTimizer, version 3.12.8), both of which are included within SCIP.

We generated a number of test problems involving $n=10,20,30$ rectangles/squares, with rectangle/square dimensions being randomly generated (to two decimal places) from $[1,5]$. For each test problem we considered three different container radii, where the container radii $R$ were set so that the area of the container
($ \pi R^2$) was 
 approximately $\frac{1}{3}$, $\frac{1}{2}$ and $\frac{2}{3}$ of the 
total area ($\sum_{i=1}^{n} L_iW_i $) of the
$n$ rectangles/squares. 
All of the randomly generated  test problems considered in this paper are publicly available from OR-Library~\cite{Beasley90}, see 
\newline
http://people.brunel.ac.uk/${\sim}$mastjjb/jeb/orlib/rspackinfo.html.

\subsection{Rectangle packing, no rotation}
Table~\ref{table2} shows the results obtained for the rectangle packing test problems considered (where the rectangles have fixed orientation, so no rotation is allowed).
In that table we show the value of $n$ and the value of the container area 
fraction. For the two objectives considered (maximise the number of rectangles packed, maximise the total area of the rectangles packed) we give the value of the best solution achieved. We also show the replication at which we first encountered the best solution shown, as well as the total time (in seconds) over all five replications.

In Table~\ref{table2}  for a fixed $n$ (and so a fixed set of rectangles to be packed) we can see that, as we would expect, as the container area fraction increases (so the container is of larger radius and we  can hence pack more of the rectangles) the solution value also increases.

As an illustration of the results obtained Figure~\ref{fig3} and Figure~\ref{fig4} show the  solutions in 
Table~\ref{table2} for the two  problems in that table with $n=30$ and the largest container area fraction.

In Figure~\ref{fig3} we can see that the solution consists of rectangles 1-18, together with rectangle 23. Recalling that rectangles are ordered in increasing size (area) order this packing is as we would expect, in that many of the smaller rectangles are used in a solution that aims to maximise the number of rectangles packed.

In Figure~\ref{fig4} we can see that the solution consists of a mix of rectangles. The first six smallest rectangles, rectangles 1-6, together with rectangles 8,11-13,17,19,22-24,28. This figure contains 16 rectangles in total, compared with the 19 rectangles used in Figure~\ref{fig3}.

\begin{table}[!htb]
{\footnotesize
{\renewcommand{\arraystretch}{1.5}
\setlength{\tabcolsep}{8pt}
\begin{tabular}{ccccccccc}
\hline
Number & Container & \multicolumn{3}{c}{Maximise number} & & \multicolumn{3}{c}{Maximise area} \\
\cline{3-5} \cline{7-9} 
of & area & Best & Replication & Total & 
& Best & Replication & Total \\
rectangles ($n$) & fraction & solution  & & time (s) & & solution & & time (s) \\
 
 \hline

10	&	$\frac{1}{3}$	&	5	&	2	&	3058	&	&	18.4441	&	1	&	3292	\\	   
	&	$\frac{1}{2}$	&	6	&	1	&	2862	&	&	28.9390	&	1	&	2992	\\	   
	&	$\frac{2}{3}$	&	7	&	1	&	2966	&	&	37.6878	&	2	&	4754	\\	   
20	&	$\frac{1}{3}$	&	7	&	1	&	6278	&	&	43.3885	&	1	&	7227	\\	   
	&	$\frac{1}{2}$	&	10	&	5	&	4530	&	&	63.1643	&	1	&	9791	\\	   
	&	$\frac{2}{3}$	&	11	&	1	&	7311	&	&	84.4446	&	2	&	10601	\\	   
30	&	$\frac{1}{3}$	&	13	&	5	&	11514	&	&	60.3570	&	4	&	14011	\\	   
	&	$\frac{1}{2}$	&	16	&	5	&	10029	&	&	85.2113	&	5	&	19786	\\	   
	&	$\frac{2}{3}$	&	19	&	5	&	6966	&	&	103.4802	&	5	&	19470	\\	 

\hline	 
											 
\end{tabular}
}}
 \caption {Computational results: rectangle packing, no rotation}\label{table2}
\end{table}

\subsection{Square packing}
Table~\ref{table3} shows the results obtained for the square packing test problems  considered. This table has the same format as Table~\ref{table2}.
As an illustration of the results obtained Figure~\ref{fig5} and Figure~\ref{fig6} show the solutions in 
Table~\ref{table3} for the two  problems in that table with $n=30$ and the largest container area fraction.

For Figure~\ref{fig5}, since we are maximising the number of squares packed, the solution must consist of the first $K$ squares, for some $K$ (the squares being ordered in increasing size order). 
In Figure~\ref{fig5} we can see that all squares up to and including  square 
$K=23$ are packed. Visually whether square 24, which must be at least as large as square 23 (and possibly larger), can also be packed into the circular container through judicious rearrangement of all of the currently positioned squares is unclear. 

In Figure~\ref{fig6}  we can see that the packing consists of a mix of squares. Squares 1-19, which are the 19 smallest squares, are all packed along with the two of the larger squares, squares 22 and 28.  This figure contains 21 squares in total, compared with the 23 squares used in Figure~\ref{fig5}.

\begin{table}[!htb]
{\footnotesize
{\renewcommand{\arraystretch}{1.5}
\setlength{\tabcolsep}{8pt}
\begin{tabular}{ccccccccc}
\hline
Number & Container & \multicolumn{3}{c}{Maximise number} & & \multicolumn{3}{c}{Maximise area} \\
\cline{3-5} \cline{7-9} 
of & area & Best & Replication & Total & 
& Best & Replication & Total \\
rectangles  ($n$) & fraction & solution  & & time (s) & & solution & & time (s) \\
 
 \hline

10	&	$\frac{1}{3}$	&	4	&	1	&	1123	&	&	22.9485	&	1	&	2762	\\	   
	&	$\frac{1}{2}$	&	5	&	1	&	2761	&	&	36.7126	&	1	&	3402	\\	   
	&	$\frac{2}{3}$	&	6	&	1	&	2275	&	&	51.7583	&	3	&	4593	\\	   
20	&	$\frac{1}{3}$	&	11	&	5	&	5450	&	&	54.1054	&	5	&	9412	\\	   
	&	$\frac{1}{2}$	&	12	&	1	&	6465	&	&	85.2107	&	4	&	11304	\\	   
	&	$\frac{2}{3}$	&	14	&	1	&	6995	&	&	109.8363	&	5	&	7636	\\	   
30	&	$\frac{1}{3}$	&	16	&	2	&	13552	&	&	54.4941	&	5	&	16629	\\	   
	&	$\frac{1}{2}$	&	20	&	2	&	13457	&	&	77.5814	&	4	&	14808	\\	   
	&	$\frac{2}{3}$	&	23	&	5	&	10427	&	&	103.0963	&	5	&	15145	\\	 
\hline	 
											 
\end{tabular}
}}
 \caption {Computational results: square packing}\label{table3}
\end{table}

\subsection{Rectangle packing, rotation allowed}
Table~\ref{table2a} shows the results obtained for the rectangle packing test problems considered in Table~\ref{table2}, but where rotation through ninety degrees is allowed. That table has the same format as Table~\ref{table2}.

As an illustration of the results obtained Figure~\ref{fig3a} and Figure~\ref{fig4a} show the  solutions in 
Table~\ref{table2} for the  two problems in that table with 30 rectangles and the largest container area fraction. In those figures the letter r after the rectangle number indicates that the rectangle has been rotated through ninety degrees.

Comparing Table~\ref{table2a} with Table~\ref{table2} we can see that the solution value where rotation is allowed is greater than (or equal to) the solution value with no rotation for all but two of the 18 test problems considered.

In the discussion above as to how to extend our formulation 
to deal with rotation through ninety degrees we noted the increase in the consequent size of the formulation, both with respect to the number of linear and nonlinear  constraints and with respect to the number of variables. Comparing the computation times in 
Table~\ref{table2a}
with those in 
Table~\ref{table2} 
 does indeed indicate 
that dealing with a problem where rectangles  can be rotated is  much more challenging computationally than dealing with a problem where the rectangles have fixed orientation.

\begin{table}[!htb]
{\footnotesize
{\renewcommand{\arraystretch}{1.5}
\setlength{\tabcolsep}{8pt}
\begin{tabular}{ccccccccc}
\hline
Number & Container & \multicolumn{3}{c}{Maximise number} & & \multicolumn{3}{c}{Maximise area} \\
\cline{3-5} \cline{7-9} 
of & area & Best & Replication & Total & 
& Best & Replication & Total \\
rectangles  & fraction & solution  & & time (s) & & solution & & time (s) \\
 
 \hline

10	&	$\frac{1}{3}$	&	5	&	1	&	9836	&	&	19.6702	&	1	&	8771	\\	   
	&	$\frac{1}{2}$	&	6	&	1	&	10332	&	&	29.5041	&	1	&	16093	\\	   
	&	$\frac{2}{3}$	&	7	&	1	&	12409	&	&	37.9687	&	2	&	15526	\\	   
20	&	$\frac{1}{3}$	&	8	&	3	&	22759	&	&	43.6850	&	2	&	50558	\\	   
	&	$\frac{1}{2}$	&	10	&	1	&	30682	&	&	63.5279	&	1	&	50013	\\	   
	&	$\frac{2}{3}$	&	12	&	4	&	30823	&	&	84.7008	&	3	&	63350	\\	   
30	&	$\frac{1}{3}$	&	14	&	1	&	49724	&	&	57.9328	&	5	&	69565	\\	   
	&	$\frac{1}{2}$	&	17	&	1	&	45857	&	&	84.3715	&	1	&	82101	\\	   
	&	$\frac{2}{3}$	&	20	&	1	&	57427	&	&	110.3253	&	3	&	39564	\\

\hline	 
											 
\end{tabular}
}}
 \caption {Computational results: rectangle packing, rotation allowed}\label{table2a}
\end{table}

\subsection{Comment}
As with many heuristic algorithms presented in the literature it is difficult to draw firm conclusions as to the quality of the results obtained without knowing either the optimal solutions of the test problems solved, or the results obtained by other heuristic algorithms  by other authors on the same set of test problems. 

For the problem considered in this paper we are not aware of any appropriate publicly available test problems which could be used to provide direct insight into the quality of our heuristic.
We would stress here however that all of the test problems used in this paper  are publicly available for use by future workers to see if they can develop approaches that perform better than the formulation space heuristic presented in this paper.

Despite this lack of appropriate test problems it is possible to gain some insight into the quality of our heuristic by taking test problems associated with a slightly different (but similar) problem. This problem is the problem of packing $n$ unit squares within a circle of small (ideally minimal) radius. Here, unlike the problem we consider,  all squares must be packed  (whereas our heuristic is particularised for the case where one or more squares need not be packed).

We used our heuristic to maximise the number of unit squares packed into a circular container of known radius utilising the test problems given by  Friedman~\cite{Friedman2018}. For these  
problems~\cite{Friedman2018} gives the best solution known for the minimum radius circle within which it is possible to pack all $n$ unit squares. 
Some of these best known solutions involve arbitrary rotation (which our heuristic cannot deal with) and so we only considered problems which  did not involve rotation. Note also here that, as far as we aware, the results given in~\cite{Friedman2018} were found by varying authors using varying approaches (including, we believe, results based on human intervention). This contrasts with our results produced by a single algorithmic heuristic approach that does not involve any human intervention.

The results are shown in Table~\ref{tablef}. In that table we show the number of unit squares ($n$) and  the value of best solution (maximum number of unit squares packed) as found by our heuristic.
We also show the replication at which we first encountered the best solution shown, as well as the total time (in seconds) over all five replications.
Considering Table~\ref{tablef} we can see that for 13 of the 16  problems considered our heuristic succeeds in 
finding the best known solution by
 packing all $n$ unit squares into the given circular container.

\begin{table}[!htb]
\centering 

\begin{tabular}{cccc}
\hline
	Number of 	&	Best solution	&	Replication & Total time (s)	\\   
	 unit squares ($n$)	&	(number of unit squares packed) 	&	 & 	\\  
 \hline

1	&	1	&	1	&	0.0	\\	   
2	&	2	&	1	&	0.2	\\	   
3	&	3	&	1	&	0.3	\\	   
4	&	4	&	1	&	0.7	\\	   
5	&	5	&	1	&	3.6	\\	   	   
7	&	7	&	1	&	140.2	\\	     
9	&	9	&	1	&	49.2	\\	   
10	&	10	&	1	&	158.5	\\	   
11	&	10	&	1	&	2619.5	\\	   
12	&	12	&	1	&	120.2	\\	      
14	&	14	&	1	&	460.6	\\	     
16	&	16	&	5	&	4800.7	\\	     
18	&	18	&	2	&	3153.6	\\	     
21	&	21	&	5	&	4961.4	\\	   	   
26	&	23	&	1	&	8756.0	\\	 
30	&	27	&	1	&	18984.5	\\	 

	\hline	 
											 
\end{tabular}
\caption{Computational results: unit square packing}\label{tablef}
\end{table}

\section{Conclusions and future work} \label{Sec:Conclusions}

In this paper we have formulated the problem of packing unequal rectangles/squares into a fixed size circular container as a 
mixed-integer nonlinear program. We showed how we can eliminate a nonlinear maximisation term that arises in one of the constraints in our formulation and 
indicated the amendments that can be made to the formulation
when considering packing squares so as to maximise the number of squares packed.

We discussed how to amend our formulation to deal with the case  where unequal rectangles can be rotated through ninety degrees.
A formulation space search heuristic was presented and
computational results given for test problems involving up to 30 rectangles/squares,  with these test problems being made publicly available for future workers.

In terms of future work we plan to investigate changes to our formulation, for example by making use of  McCormick cuts to replace products of variables.

\newpage

\begin{figure}[!htb]
\centering
\begin{tikzpicture}[scale=1];
\draw (0,0) circle [radius=        4.0000000000cm];
\draw(       -3.7096403558 ,         -.4538739448) -- 
(       -2.8517016371 ,         -.4538739448) -- 
(       -2.8517016371 ,        -1.3173836941) -- 
(       -3.7096403558 ,        -1.3173836941) -- 
(       -3.7096403558 ,         -.4538739448) ;
 \node at (       -3.2806709965 ,         -.8856288195) {     1}; 
\draw(         .4237081360 ,         -.8756024043) -- 
(        1.1312290273 ,         -.8756024043) -- 
(        1.1312290273 ,        -1.9285271954) -- 
(         .4237081360 ,        -1.9285271954) -- 
(         .4237081360 ,         -.8756024043) ;
 \node at (         .7774685816 ,        -1.4020647999) {     2}; 
\draw(        1.1341336498 ,         -.7308883810) -- 
(        1.8249414492 ,         -.7308883810) -- 
(        1.8249414492 ,        -1.8729496624) -- 
(        1.1341336498 ,        -1.8729496624) -- 
(        1.1341336498 ,         -.7308883810) ;
 \node at (        1.4795375495 ,        -1.3019190217) {     3}; 
\draw(         .2334830750 ,        -2.5839516405) -- 
(        1.3476892031 ,        -2.5839516405) -- 
(        1.3476892031 ,        -3.3304697463) -- 
(         .2334830750 ,        -3.3304697463) -- 
(         .2334830750 ,        -2.5839516405) ;
 \node at (         .7905861390 ,        -2.9572106934) {     4}; 
\draw(        1.4156759167 ,          .6695458990) -- 
(        2.7638653317 ,          .6695458990) -- 
(        2.7638653317 ,          .0511614979) -- 
(        1.4156759167 ,          .0511614979) -- 
(        1.4156759167 ,          .6695458990) ;
 \node at (        2.0897706242 ,          .3603536984) {     5}; 
\draw(        1.3696208549 ,        -2.2174513363) -- 
(        2.3612643090 ,        -2.2174513363) -- 
(        2.3612643090 ,        -3.2202368516) -- 
(        1.3696208549 ,        -3.2202368516) -- 
(        1.3696208549 ,        -2.2174513363) ;
 \node at (        1.8654425820 ,        -2.7188440939) {     6}; 
\draw(       -1.8093054342 ,        -1.6201055717) -- 
(        -.2717009774 ,        -1.6201055717) -- 
(        -.2717009774 ,        -2.3276264631) -- 
(       -1.8093054342 ,        -2.3276264631) -- 
(       -1.8093054342 ,        -1.6201055717) ;
 \node at (       -1.0405032058 ,        -1.9738660174) {     7}; 
\draw(       -2.9216206652 ,         2.7251311734) -- 
(       -1.7517042307 ,         2.7251311734) -- 
(       -1.7517042307 ,         1.7112035968) -- 
(       -2.9216206652 ,         1.7112035968) -- 
(       -2.9216206652 ,         2.7251311734) ;
 \node at (       -2.3366624479 ,         2.2181673851) {     8}; 
\draw(         .9848750858 ,         -.1265125474) -- 
(        3.4194154758 ,         -.1265125474) -- 
(        3.4194154758 ,         -.6947576728) -- 
(         .9848750858 ,         -.6947576728) -- 
(         .9848750858 ,         -.1265125474) ;
 \node at (        2.2021452808 ,         -.4106351101) {     9}; 
\draw(         .5963384613 ,         1.6166460329) -- 
(        2.1339429182 ,         1.6166460329) -- 
(        2.1339429182 ,          .6974259772) -- 
(         .5963384613 ,          .6974259772) -- 
(         .5963384613 ,         1.6166460329) ;
 \node at (        1.3651406897 ,         1.1570360051) {    10}; 
\draw(        -.2697592050 ,         -.1744623059) -- 
(         .4210485945 ,         -.1744623059) -- 
(         .4210485945 ,        -2.4920110524) -- 
(        -.2697592050 ,        -2.4920110524) -- 
(        -.2697592050 ,         -.1744623059) ;
 \node at (         .0756446948 ,        -1.3332366791) {    11}; 
\draw(       -1.6521713097 ,        -2.5562339823) -- 
(        -.0477144852 ,        -2.5562339823) -- 
(        -.0477144852 ,        -3.5757325895) -- 
(       -1.6521713097 ,        -3.5757325895) -- 
(       -1.6521713097 ,        -2.5562339823) ;
 \node at (        -.8499428974 ,        -3.0659832859) {    12}; 
\draw(        -.9178102520 ,          .6952737631) -- 
(        1.1490421157 ,          .6952737631) -- 
(        1.1490421157 ,         -.1236677410) -- 
(        -.9178102520 ,         -.1236677410) -- 
(        -.9178102520 ,          .6952737631) ;
 \node at (         .1156159319 ,          .2858030111) {    13}; 
\draw(       -3.6102829286 ,          .6929869857) -- 
(       -1.0420378033 ,          .6929869857) -- 
(       -1.0420378033 ,         -.1203834879) -- 
(       -3.6102829286 ,         -.1203834879) -- 
(       -3.6102829286 ,          .6929869857) ;
 \node at (       -2.3261603660 ,          .2863017489) {    14}; 
\draw(        1.8277993625 ,         -.6983470946) -- 
(        3.3375486661 ,         -.6983470946) -- 
(        3.3375486661 ,        -2.2025253676) -- 
(        1.8277993625 ,        -2.2025253676) -- 
(        1.8277993625 ,         -.6983470946) ;
 \node at (        2.5826740143 ,        -1.4504362311) {    15}; 
\draw(       -3.6196305670 ,         1.6972157268) -- 
(        -.9566779207 ,         1.6972157268) -- 
(        -.9566779207 ,          .7000012421) -- 
(       -3.6196305670 ,          .7000012421) -- 
(       -3.6196305670 ,         1.6972157268) ;
 \node at (       -2.2881542439 ,         1.1986084845) {    16}; 
\draw(       -2.8195899443 ,         -.1229780119) -- 
(       -1.8112333983 ,         -.1229780119) -- 
(       -1.8112333983 ,        -2.8360699339) -- 
(       -2.8195899443 ,        -2.8360699339) -- 
(       -2.8195899443 ,         -.1229780119) ;
 \node at (       -2.3154116713 ,        -1.4795239729) {    17}; 
\draw(        -.9541285408 ,         3.5013216232) -- 
(         .0597990358 ,         3.5013216232) -- 
(         .0597990358 ,          .7548035174) -- 
(        -.9541285408 ,          .7548035174) -- 
(        -.9541285408 ,         3.5013216232) ;
 \node at (        -.4471647525 ,         2.1280625703) {    18}; 
\draw(         .0624519736 ,         3.3676310334) -- 
(        2.1571594945 ,         3.3676310334) -- 
(        2.1571594945 ,         1.6183274122) -- 
(         .0624519736 ,         1.6183274122) -- 
(         .0624519736 ,         3.3676310334) ;
 \node at (        1.1098057340 ,         2.4929792228) {    23}; 
\end{tikzpicture}
\captionsetup{font=small}
\caption{Maximising number of rectangles packed, no rotation, $n= 30$, container area fraction $\frac{2}{3}$ }
\label{fig3}
\end{figure}
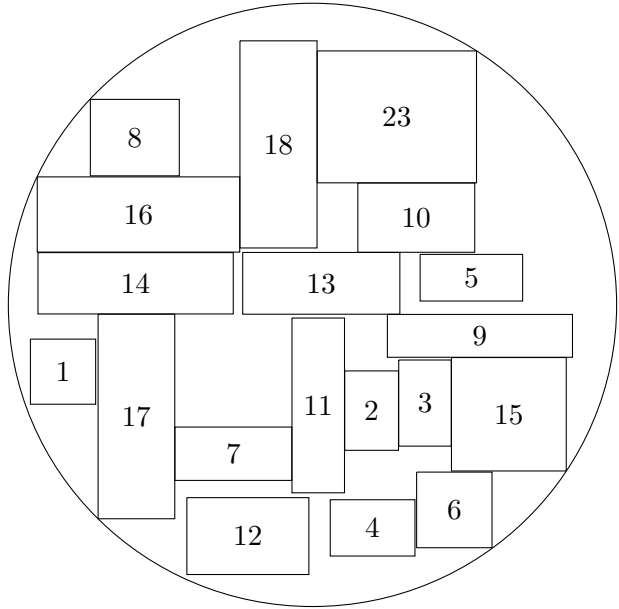

\begin{figure}[!htb]
\centering
\begin{tikzpicture}[scale=1];
\draw (0,0) circle [radius=        4.0000000000cm];
\draw(         .9088831876 ,         1.4584689091) -- 
(        1.7668219063 ,         1.4584689091) -- 
(        1.7668219063 ,          .5949591598) -- 
(         .9088831876 ,          .5949591598) -- 
(         .9088831876 ,         1.4584689091) ;
 \node at (        1.3378525469 ,         1.0267140344) {     1}; 
\draw(         .4874594449 ,         -.3298385613) -- 
(        1.1949803363 ,         -.3298385613) -- 
(        1.1949803363 ,        -1.3827633524) -- 
(         .4874594449 ,        -1.3827633524) -- 
(         .4874594449 ,         -.3298385613) ;
 \node at (         .8412198906 ,         -.8563009569) {     2}; 
\draw(        -.1387332466 ,        -1.8064702814) -- 
(         .5520745528 ,        -1.8064702814) -- 
(         .5520745528 ,        -2.9485315628) -- 
(        -.1387332466 ,        -2.9485315628) -- 
(        -.1387332466 ,        -1.8064702814) ;
 \node at (         .2066706531 ,        -2.3775009221) {     3}; 
\draw(         .8204467519 ,         3.4786869205) -- 
(        1.9346528800 ,         3.4786869205) -- 
(        1.9346528800 ,         2.7321688146) -- 
(         .8204467519 ,         2.7321688146) -- 
(         .8204467519 ,         3.4786869205) ;
 \node at (        1.3775498160 ,         3.1054278676) {     4}; 
\draw(       -2.1482918228 ,          .1167806288) -- 
(        -.8001024078 ,          .1167806288) -- 
(        -.8001024078 ,         -.5016037723) -- 
(       -2.1482918228 ,         -.5016037723) -- 
(       -2.1482918228 ,          .1167806288) ;
 \node at (       -1.4741971153 ,         -.1924115718) {     5}; 
\draw(       -2.0690504767 ,         -.5258576968) -- 
(       -1.0774070227 ,         -.5258576968) -- 
(       -1.0774070227 ,        -1.5286432121) -- 
(       -2.0690504767 ,        -1.5286432121) -- 
(       -2.0690504767 ,         -.5258576968) ;
 \node at (       -1.5732287497 ,        -1.0272504544) {     6}; 
\draw(       -1.5577106395 ,         3.5802664894) -- 
(        -.3877942049 ,         3.5802664894) -- 
(        -.3877942049 ,         2.5663389128) -- 
(       -1.5577106395 ,         2.5663389128) -- 
(       -1.5577106395 ,         3.5802664894) ;
 \node at (        -.9727524222 ,         3.0733027011) {     8}; 
\draw(         .2148130976 ,         2.4975667791) -- 
(         .9056208970 ,         2.4975667791) -- 
(         .9056208970 ,          .1800180326) -- 
(         .2148130976 ,          .1800180326) -- 
(         .2148130976 ,         2.4975667791) ;
 \node at (         .5602169973 ,         1.3387924059) {    11}; 
\draw(        1.1621506497 ,         2.6819252186) -- 
(        2.7666074743 ,         2.6819252186) -- 
(        2.7666074743 ,         1.6624266114) -- 
(        1.1621506497 ,         1.6624266114) -- 
(        1.1621506497 ,         2.6819252186) ;
 \node at (        1.9643790620 ,         2.1721759150) {    12}; 
\draw(        1.2915022243 ,         -.7867065239) -- 
(        3.3583545919 ,         -.7867065239) -- 
(        3.3583545919 ,        -1.6056480281) -- 
(        1.2915022243 ,        -1.6056480281) -- 
(        1.2915022243 ,         -.7867065239) ;
 \node at (        2.3249284081 ,        -1.1961772760) {    13}; 
\draw(        -.7968619713 ,         1.2032321938) -- 
(         .2114945747 ,         1.2032321938) -- 
(         .2114945747 ,        -1.5098597282) -- 
(        -.7968619713 ,        -1.5098597282) -- 
(        -.7968619713 ,         1.2032321938) ;
 \node at (        -.2926836983 ,         -.1533137672) {    17}; 
\draw(         .5641138320 ,        -1.6138444828) -- 
(        2.5808269239 ,        -1.6138444828) -- 
(        2.5808269239 ,        -3.0511703881) -- 
(         .5641138320 ,        -3.0511703881) -- 
(         .5641138320 ,        -1.6138444828) ;
 \node at (        1.5724703779 ,        -2.3325074354) {    19}; 
\draw(       -3.7671460247 ,          .8005831732) -- 
(       -2.1515471389 ,          .8005831732) -- 
(       -2.1515471389 ,        -1.3386925928) -- 
(       -3.7671460247 ,        -1.3386925928) -- 
(       -3.7671460247 ,          .8005831732) ;
 \node at (       -2.9593465818 ,         -.2690547098) {    22}; 
\draw(       -2.2523767364 ,        -1.5498963437) -- 
(        -.1576692155 ,        -1.5498963437) -- 
(        -.1576692155 ,        -3.2991999649) -- 
(       -2.2523767364 ,        -3.2991999649) -- 
(       -2.2523767364 ,        -1.5498963437) ;
 \node at (       -1.2050229760 ,        -2.4245481543) {    23}; 
\draw(       -3.0558099829 ,         2.5416326818) -- 
(        -.8496818492 ,         2.5416326818) -- 
(        -.8496818492 ,          .8201842138) -- 
(       -3.0558099829 ,          .8201842138) -- 
(       -3.0558099829 ,         2.5416326818) ;
 \node at (       -1.9527459161 ,         1.6809084478) {    24}; 
\draw(        1.7702059005 ,         1.6120274523) -- 
(        3.6587852877 ,         1.6120274523) -- 
(        3.6587852877 ,         -.7779446926) -- 
(        1.7702059005 ,         -.7779446926) -- 
(        1.7702059005 ,         1.6120274523) ;
 \node at (        2.7144955941 ,          .4170413799) {    28}; 
\end{tikzpicture}
\captionsetup{font=small}
\caption{Maximising total area of rectangles packed, no rotation,  $n= 30$, container area fraction $\frac{2}{3}$ }
\label{fig4}
\end{figure}
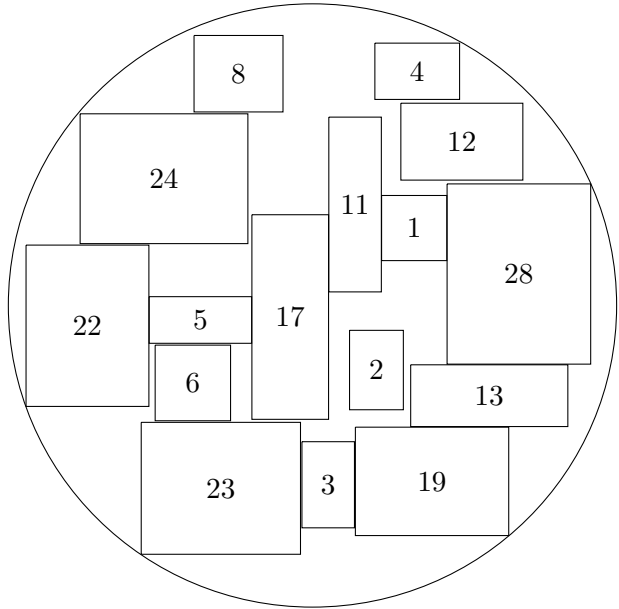

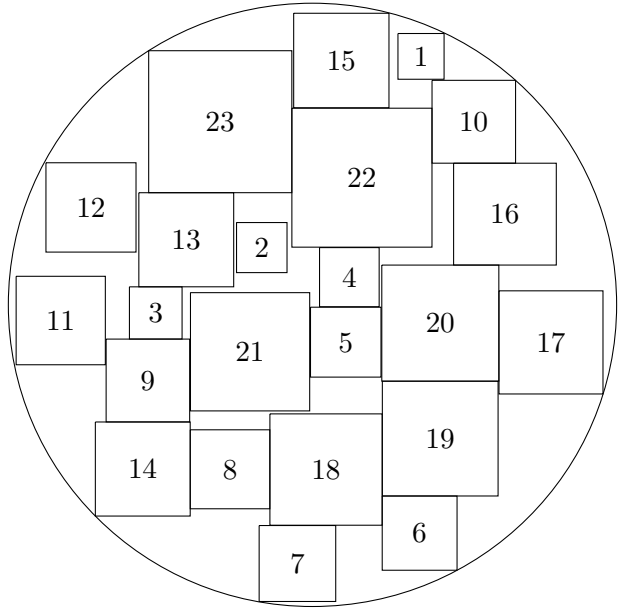
\begin{figure}[!htb]
\centering
\begin{tikzpicture}[scale=1];
\draw (0,0) circle [radius=        4.0000000000cm];
\draw(        1.1240872451 ,         3.5982473001) -- 
(        1.7306680748 ,         3.5982473001) -- 
(        1.7306680748 ,         2.9916664704) -- 
(        1.1240872451 ,         2.9916664704) -- 
(        1.1240872451 ,         3.5982473001) ;
 \node at (        1.4273776599 ,         3.2949568852) {     1}; 
\draw(        -.9987322003 ,         1.0913705409) -- 
(        -.3349267640 ,         1.0913705409) -- 
(        -.3349267640 ,          .4275651046) -- 
(        -.9987322003 ,          .4275651046) -- 
(        -.9987322003 ,         1.0913705409) ;
 \node at (        -.6668294822 ,          .7594678228) {     2}; 
\draw(       -2.4075410268 ,          .2389884344) -- 
(       -1.7151232872 ,          .2389884344) -- 
(       -1.7151232872 ,         -.4534293052) -- 
(       -2.4075410268 ,         -.4534293052) -- 
(       -2.4075410268 ,          .2389884344) ;
 \node at (       -2.0613321570 ,         -.1072204354) {     3}; 
\draw(         .0928802237 ,          .7592135089) -- 
(         .8768573338 ,          .7592135089) -- 
(         .8768573338 ,         -.0247636012) -- 
(         .0928802237 ,         -.0247636012) -- 
(         .0928802237 ,          .7592135089) ;
 \node at (         .4848687787 ,          .3672249538) {     4}; 
\draw(        -.0260879301 ,         -.0310581519) -- 
(         .9009506965 ,         -.0310581519) -- 
(         .9009506965 ,         -.9580967785) -- 
(        -.0260879301 ,         -.9580967785) -- 
(        -.0260879301 ,         -.0310581519) ;
 \node at (         .4374313832 ,         -.4945774652) {     5}; 
\draw(         .9170567289 ,        -2.5347852949) -- 
(        1.9013199621 ,        -2.5347852949) -- 
(        1.9013199621 ,        -3.5190485281) -- 
(         .9170567289 ,        -3.5190485281) -- 
(         .9170567289 ,        -2.5347852949) ;
 \node at (        1.4091883455 ,        -3.0269169115) {     6}; 
\draw(        -.7014888347 ,        -2.9262494663) -- 
(         .3056642411 ,        -2.9262494663) -- 
(         .3056642411 ,        -3.9334025421) -- 
(        -.7014888347 ,        -3.9334025421) -- 
(        -.7014888347 ,        -2.9262494663) ;
 \node at (        -.1979122968 ,        -3.4298260042) {     7}; 
\draw(       -1.6077258591 ,        -1.6568080075) -- 
(        -.5605155587 ,        -1.6568080075) -- 
(        -.5605155587 ,        -2.7040183079) -- 
(       -1.6077258591 ,        -2.7040183079) -- 
(       -1.6077258591 ,        -1.6568080075) ;
 \node at (       -1.0841207089 ,        -2.1804131577) {     8}; 
\draw(       -2.7139556888 ,         -.4539034160) -- 
(       -1.6152432424 ,         -.4539034160) -- 
(       -1.6152432424 ,        -1.5526158623) -- 
(       -2.7139556888 ,        -1.5526158623) -- 
(       -2.7139556888 ,         -.4539034160) ;
 \node at (       -2.1645994656 ,        -1.0032596391) {     9}; 
\draw(        1.5721596421 ,         2.9774040252) -- 
(        2.6708720884 ,         2.9774040252) -- 
(        2.6708720884 ,         1.8786915788) -- 
(        1.5721596421 ,         1.8786915788) -- 
(        1.5721596421 ,         2.9774040252) ;
 \node at (        2.1215158652 ,         2.4280478020) {    10}; 
\draw(       -3.8974526949 ,          .3787977691) -- 
(       -2.7243482600 ,          .3787977691) -- 
(       -2.7243482600 ,         -.7943066658) -- 
(       -3.8974526949 ,         -.7943066658) -- 
(       -3.8974526949 ,          .3787977691) ;
 \node at (       -3.3109004775 ,         -.2077544483) {    11}; 
\draw(       -3.5053765188 ,         1.8843399161) -- 
(       -2.3208271626 ,         1.8843399161) -- 
(       -2.3208271626 ,          .6997905599) -- 
(       -3.5053765188 ,          .6997905599) -- 
(       -3.5053765188 ,         1.8843399161) ;
 \node at (       -2.9131018407 ,         1.2920652380) {    12}; 
\draw(       -2.2836941851 ,         1.4869636662) -- 
(       -1.0361977616 ,         1.4869636662) -- 
(       -1.0361977616 ,          .2394672427) -- 
(       -2.2836941851 ,          .2394672427) -- 
(       -2.2836941851 ,         1.4869636662) ;
 \node at (       -1.6599459734 ,          .8632154544) {    13}; 
\draw(       -2.8557399830 ,        -1.5530890834) -- 
(       -1.6082435595 ,        -1.5530890834) -- 
(       -1.6082435595 ,        -2.8005855069) -- 
(       -2.8557399830 ,        -2.8005855069) -- 
(       -2.8557399830 ,        -1.5530890834) ;
 \node at (       -2.2319917713 ,        -2.1768372952) {    14}; 
\draw(        -.2473599232 ,         3.8667163343) -- 
(        1.0058589610 ,         3.8667163343) -- 
(        1.0058589610 ,         2.6134974501) -- 
(        -.2473599232 ,         2.6134974501) -- 
(        -.2473599232 ,         3.8667163343) ;
 \node at (         .3792495189 ,         3.2401068922) {    15}; 
\draw(        1.8559504373 ,         1.8783243604) -- 
(        3.2064511527 ,         1.8783243604) -- 
(        3.2064511527 ,          .5278236451) -- 
(        1.8559504373 ,          .5278236451) -- 
(        1.8559504373 ,         1.8783243604) ;
 \node at (        2.5312007950 ,         1.2030740028) {    16}; 
\draw(        2.4528795813 ,          .1860678199) -- 
(        3.8205476786 ,          .1860678199) -- 
(        3.8205476786 ,        -1.1816002774) -- 
(        2.4528795813 ,        -1.1816002774) -- 
(        2.4528795813 ,          .1860678199) ;
 \node at (        3.1367136300 ,         -.4977662287) {    17}; 
\draw(        -.5599979747 ,        -1.4456810157) -- 
(         .9163968751 ,        -1.4456810157) -- 
(         .9163968751 ,        -2.9220758654) -- 
(        -.5599979747 ,        -2.9220758654) -- 
(        -.5599979747 ,        -1.4456810157) ;
 \node at (         .1781994502 ,        -2.1838784405) {    18}; 
\draw(         .9186998882 ,        -1.0122484334) -- 
(        2.4408744233 ,        -1.0122484334) -- 
(        2.4408744233 ,        -2.5344229685) -- 
(         .9186998882 ,        -2.5344229685) -- 
(         .9186998882 ,        -1.0122484334) ;
 \node at (        1.6797871557 ,        -1.7733357009) {    19}; 
\draw(         .9114939769 ,          .5274566178) -- 
(        2.4508358939 ,          .5274566178) -- 
(        2.4508358939 ,        -1.0118852992) -- 
(         .9114939769 ,        -1.0118852992) -- 
(         .9114939769 ,          .5274566178) ;
 \node at (        1.6811649354 ,         -.2422143407) {    20}; 
\draw(       -1.6046812847 ,          .1617622537) -- 
(        -.0367270644 ,          .1617622537) -- 
(        -.0367270644 ,        -1.4061919666) -- 
(       -1.6046812847 ,        -1.4061919666) -- 
(       -1.6046812847 ,          .1617622537) ;
 \node at (        -.8207041745 ,         -.6222148564) {    21}; 
\draw(        -.2710381259 ,         2.6077724767) -- 
(        1.5715942060 ,         2.6077724767) -- 
(        1.5715942060 ,          .7651401448) -- 
(        -.2710381259 ,          .7651401448) -- 
(        -.2710381259 ,         2.6077724767) ;
 \node at (         .6502780400 ,         1.6864563107) {    22}; 
\draw(       -2.1543038971 ,         3.3701279201) -- 
(        -.2716143406 ,         3.3701279201) -- 
(        -.2716143406 ,         1.4874383636) -- 
(       -2.1543038971 ,         1.4874383636) -- 
(       -2.1543038971 ,         3.3701279201) ;
 \node at (       -1.2129591188 ,         2.4287831418) {    23}; 
\end{tikzpicture}
\caption{Maximising number of squares packed, $n= 30$, container area fraction $\frac{2}{3}$ }
\label{fig5}
\end{figure}

\begin{figure}[!htb]
\centering
\begin{tikzpicture}[scale=1];
\draw (0,0) circle [radius=        4.0000000000cm];
\draw(       -2.2686252086 ,         3.2802816032) -- 
(       -1.6620443789 ,         3.2802816032) -- 
(       -1.6620443789 ,         2.6737007735) -- 
(       -2.2686252086 ,         2.6737007735) -- 
(       -2.2686252086 ,         3.2802816032) ;
 \node at (       -1.9653347937 ,         2.9769911884) {     1}; 
\draw(        -.5708802493 ,         1.5321560821) -- 
(         .0929251871 ,         1.5321560821) -- 
(         .0929251871 ,          .8683506458) -- 
(        -.5708802493 ,          .8683506458) -- 
(        -.5708802493 ,         1.5321560821) ;
 \node at (        -.2389775311 ,         1.2002533639) {     2}; 
\draw(        -.4220350272 ,         3.2944574933) -- 
(         .2703827124 ,         3.2944574933) -- 
(         .2703827124 ,         2.6020397537) -- 
(        -.4220350272 ,         2.6020397537) -- 
(        -.4220350272 ,         3.2944574933) ;
 \node at (        -.0758261574 ,         2.9482486235) {     3}; 
\draw(       -1.4112045662 ,         1.9914838227) -- 
(        -.6272274560 ,         1.9914838227) -- 
(        -.6272274560 ,         1.2075067125) -- 
(       -1.4112045662 ,         1.2075067125) -- 
(       -1.4112045662 ,         1.9914838227) ;
 \node at (       -1.0192160111 ,         1.5994952676) {     4}; 
\draw(       -3.5895689525 ,         1.7633945342) -- 
(       -2.6625303259 ,         1.7633945342) -- 
(       -2.6625303259 ,          .8363559076) -- 
(       -3.5895689525 ,          .8363559076) -- 
(       -3.5895689525 ,         1.7633945342) ;
 \node at (       -3.1260496392 ,         1.2998752209) {     5}; 
\draw(       -1.0916952483 ,          .8602959263) -- 
(        -.1074320151 ,          .8602959263) -- 
(        -.1074320151 ,         -.1239673069) -- 
(       -1.0916952483 ,         -.1239673069) -- 
(       -1.0916952483 ,          .8602959263) ;
 \node at (        -.5995636317 ,          .3681643097) {     6}; 
\draw(        -.4372634338 ,         2.5477549264) -- 
(         .5698896420 ,         2.5477549264) -- 
(         .5698896420 ,         1.5406018506) -- 
(        -.4372634338 ,         1.5406018506) -- 
(        -.4372634338 ,         2.5477549264) ;
 \node at (         .0663131041 ,         2.0441783885) {     7}; 
\draw(       -2.6610445314 ,         2.6248359795) -- 
(       -1.6138342310 ,         2.6248359795) -- 
(       -1.6138342310 ,         1.5776256790) -- 
(       -2.6610445314 ,         1.5776256790) -- 
(       -2.6610445314 ,         2.6248359795) ;
 \node at (       -2.1374393812 ,         2.1012308293) {     8}; 
\draw(        2.8519405801 ,          .4794214330) -- 
(        3.9506530264 ,          .4794214330) -- 
(        3.9506530264 ,         -.6192910133) -- 
(        2.8519405801 ,         -.6192910133) -- 
(        2.8519405801 ,          .4794214330) ;
 \node at (        3.4012968032 ,         -.0699347901) {     9}; 
\draw(        1.9409375217 ,         2.5991089089) -- 
(        3.0396499680 ,         2.5991089089) -- 
(        3.0396499680 ,         1.5003964625) -- 
(        1.9409375217 ,         1.5003964625) -- 
(        1.9409375217 ,         2.5991089089) ;
 \node at (        2.4902937449 ,         2.0497526857) {    10}; 
\draw(       -1.6123233497 ,         3.3179468180) -- 
(        -.4392189148 ,         3.3179468180) -- 
(        -.4392189148 ,         2.1448423831) -- 
(       -1.6123233497 ,         2.1448423831) -- 
(       -1.6123233497 ,         3.3179468180) ;
 \node at (       -1.0257711323 ,         2.7313946006) {    11}; 
\draw(       -3.7825010330 ,          .8332563403) -- 
(       -2.5979516768 ,          .8332563403) -- 
(       -2.5979516768 ,         -.3512930160) -- 
(       -3.7825010330 ,         -.3512930160) -- 
(       -3.7825010330 ,          .8332563403) ;
 \node at (       -3.1902263549 ,          .2409816621) {    12}; 
\draw(       -2.4160072315 ,         -.3403811696) -- 
(       -1.1685108081 ,         -.3403811696) -- 
(       -1.1685108081 ,        -1.5878775931) -- 
(       -2.4160072315 ,        -1.5878775931) -- 
(       -2.4160072315 ,         -.3403811696) ;
 \node at (       -1.7922590198 ,         -.9641293814) {    13}; 
\draw(       -3.6646236130 ,         -.3544730669) -- 
(       -2.4171271895 ,         -.3544730669) -- 
(       -2.4171271895 ,        -1.6019694904) -- 
(       -3.6646236130 ,        -1.6019694904) -- 
(       -3.6646236130 ,         -.3544730669) ;
 \node at (       -3.0408754012 ,         -.9782212787) {    14}; 
\draw(        -.0800679186 ,        -2.5695771892) -- 
(        1.1731509655 ,        -2.5695771892) -- 
(        1.1731509655 ,        -3.8227960734) -- 
(        -.0800679186 ,        -3.8227960734) -- 
(        -.0800679186 ,        -2.5695771892) ;
 \node at (         .5465415234 ,        -3.1961866313) {    15}; 
\draw(       -2.5549389999 ,         1.0359434382) -- 
(       -1.2044382846 ,         1.0359434382) -- 
(       -1.2044382846 ,         -.3145572771) -- 
(       -2.5549389999 ,         -.3145572771) -- 
(       -2.5549389999 ,         1.0359434382) ;
 \node at (       -1.8796886422 ,          .3606930805) {    16}; 
\draw(         .5715846835 ,         3.1207760762) -- 
(        1.9392527808 ,         3.1207760762) -- 
(        1.9392527808 ,         1.7531079789) -- 
(         .5715846835 ,         1.7531079789) -- 
(         .5715846835 ,         3.1207760762) ;
 \node at (        1.2554187322 ,         2.4369420276) {    17}; 
\draw(       -1.1674036069 ,         -.1463676618) -- 
(         .3089912429 ,         -.1463676618) -- 
(         .3089912429 ,        -1.6227625116) -- 
(       -1.1674036069 ,        -1.6227625116) -- 
(       -1.1674036069 ,         -.1463676618) ;
 \node at (        -.4292061820 ,         -.8845650867) {    18}; 
\draw(        1.0763940666 ,         1.4980540617) -- 
(        2.5985686016 ,         1.4980540617) -- 
(        2.5985686016 ,         -.0241204733) -- 
(        1.0763940666 ,         -.0241204733) -- 
(        1.0763940666 ,         1.4980540617) ;
 \node at (        1.8374813341 ,          .7369667942) {    19}; 
\draw(       -1.9318017849 ,        -1.6478489806) -- 
(        -.0891694530 ,        -1.6478489806) -- 
(        -.0891694530 ,        -3.4904813125) -- 
(       -1.9318017849 ,        -3.4904813125) -- 
(       -1.9318017849 ,        -1.6478489806) ;
 \node at (       -1.0104856189 ,        -2.5691651465) {    22}; 
\draw(         .3100788768 ,         -.0264747193) -- 
(        2.8508514090 ,         -.0264747193) -- 
(        2.8508514090 ,        -2.5672472515) -- 
(         .3100788768 ,        -2.5672472515) -- 
(         .3100788768 ,         -.0264747193) ;
 \node at (        1.5804651429 ,        -1.2968609854) {    28}; 
\end{tikzpicture}
\caption{Maximising total area of squares packed, $n= 30$, container area fraction $\frac{2}{3}$ }
\label{fig6}
\end{figure}
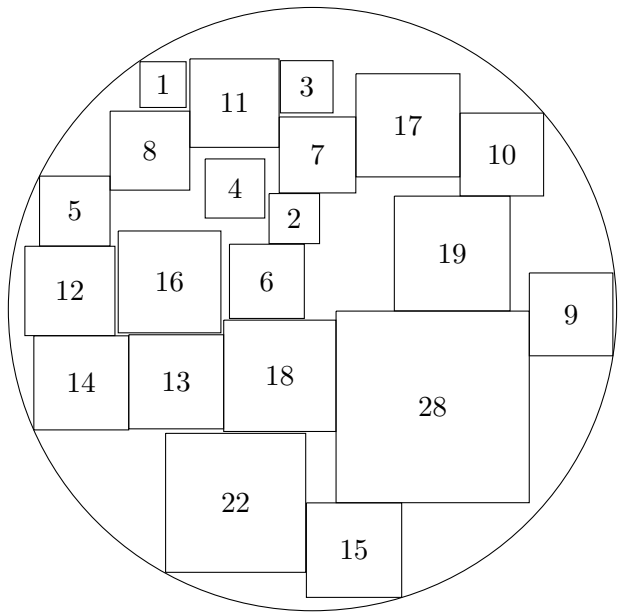

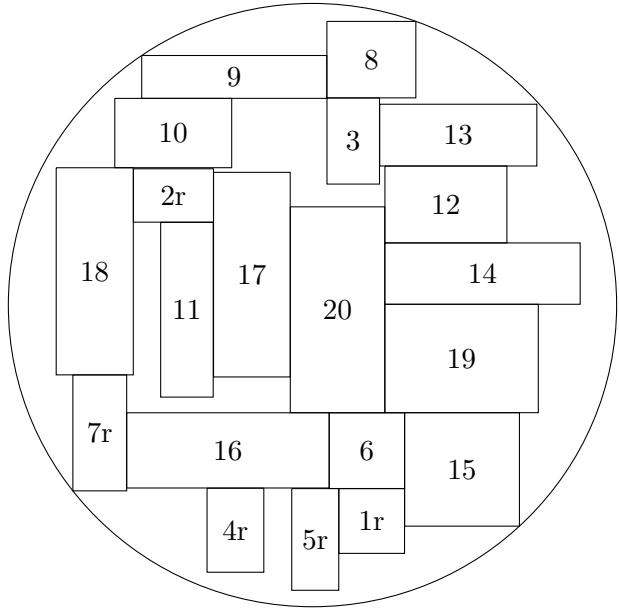
\begin{figure}[!htb]
\centering
\begin{tikzpicture}[scale=1];
\draw (0,0) circle [radius=        4.0000000000cm];
\draw(         .1906390049 ,         2.7467864217) -- 
(         .8814468044 ,         2.7467864217) -- 
(         .8814468044 ,         1.6047251403) -- 
(         .1906390049 ,         1.6047251403) -- 
(         .1906390049 ,         2.7467864217) ;
 \node at (         .5360429047 ,         2.1757557810) {     3     }; 
\draw(         .2199199479 ,        -1.4292973034) -- 
(        1.2115634019 ,        -1.4292973034) -- 
(        1.2115634019 ,        -2.4320828187) -- 
(         .2199199479 ,        -2.4320828187) -- 
(         .2199199479 ,        -1.4292973034) ;
 \node at (         .7157416749 ,        -1.9306900610) {     6     }; 
\draw(         .1894056969 ,         3.7617943428) -- 
(        1.3593221314 ,         3.7617943428) -- 
(        1.3593221314 ,         2.7478667662) -- 
(         .1894056969 ,         2.7478667662) -- 
(         .1894056969 ,         3.7617943428) ;
 \node at (         .7743639142 ,         3.2548305545) {     8     }; 
\draw(       -2.2460346569 ,         3.3097056244) -- 
(         .1885057331 ,         3.3097056244) -- 
(         .1885057331 ,         2.7414604990) -- 
(       -2.2460346569 ,         2.7414604990) -- 
(       -2.2460346569 ,         3.3097056244) ;
 \node at (       -1.0287644619 ,         3.0255830617) {     9     }; 
\draw(       -2.5999859631 ,         2.7410148394) -- 
(       -1.0623815062 ,         2.7410148394) -- 
(       -1.0623815062 ,         1.8217947837) -- 
(       -2.5999859631 ,         1.8217947837) -- 
(       -2.5999859631 ,         2.7410148394) ;
 \node at (       -1.8311837346 ,         2.2814048115) {    10    }; 
\draw(       -1.9965866279 ,         1.0974240742) -- 
(       -1.3057788284 ,         1.0974240742) -- 
(       -1.3057788284 ,        -1.2201246723) -- 
(       -1.9965866279 ,        -1.2201246723) -- 
(       -1.9965866279 ,         1.0974240742) ;
 \node at (       -1.6511827282 ,         -.0613502991) {    11    }; 
\draw(         .9526169807 ,         1.8441567633) -- 
(        2.5570738052 ,         1.8441567633) -- 
(        2.5570738052 ,          .8246581560) -- 
(         .9526169807 ,          .8246581560) -- 
(         .9526169807 ,         1.8441567633) ;
 \node at (        1.7548453930 ,         1.3344074596) {    12   }; 
\draw(         .8838049572 ,         2.6634622769) -- 
(        2.9506573249 ,         2.6634622769) -- 
(        2.9506573249 ,         1.8445207727) -- 
(         .8838049572 ,         1.8445207727) -- 
(         .8838049572 ,         2.6634622769) ;
 \node at (        1.9172311411 ,         2.2539915248) {    13    }; 
\draw(         .9523179121 ,          .8243351445) -- 
(        3.5205630375 ,          .8243351445) -- 
(        3.5205630375 ,          .0109646710) -- 
(         .9523179121 ,          .0109646710) -- 
(         .9523179121 ,          .8243351445) ;
 \node at (        2.2364404748 ,          .4176499078) {    14    }; 
\draw(        1.2119692034 ,        -1.4269515126) -- 
(        2.7217185070 ,        -1.4269515126) -- 
(        2.7217185070 ,        -2.9311297856) -- 
(        1.2119692034 ,        -2.9311297856) -- 
(        1.2119692034 ,        -1.4269515126) ;
 \node at (        1.9668438552 ,        -2.1790406491) {    15    }; 
\draw(       -2.4434279873 ,        -1.4280680143) -- 
(         .2195246589 ,        -1.4280680143) -- 
(         .2195246589 ,        -2.4252824990) -- 
(       -2.4434279873 ,        -2.4252824990) -- 
(       -2.4434279873 ,        -1.4280680143) ;
 \node at (       -1.1119516642 ,        -1.9266752566) {    16    }; 
\draw(       -1.3011216585 ,         1.7603180301) -- 
(        -.2927651125 ,         1.7603180301) -- 
(        -.2927651125 ,         -.9527738919) -- 
(       -1.3011216585 ,         -.9527738919) -- 
(       -1.3011216585 ,         1.7603180301) ;
 \node at (        -.7969433855 ,          .4037720691) {    17    }; 
\draw(       -3.3705357935 ,         1.8212656096) -- 
(       -2.3566082169 ,         1.8212656096) -- 
(       -2.3566082169 ,         -.9252524963) -- 
(       -3.3705357935 ,         -.9252524963) -- 
(       -3.3705357935 ,         1.8212656096) ;
 \node at (       -2.8635720052 ,          .4480065567) {    18    }; 
\draw(         .9530443117 ,          .0106703559) -- 
(        2.9697574037 ,          .0106703559) -- 
(        2.9697574037 ,        -1.4266555493) -- 
(         .9530443117 ,        -1.4266555493) -- 
(         .9530443117 ,          .0106703559) ;
 \node at (        1.9614008577 ,         -.7079925967) {    19    }; 
\draw(        -.2919069923 ,         1.3029313232) -- 
(         .9504328406 ,         1.3029313232) -- 
(         .9504328406 ,        -1.4268736907) -- 
(        -.2919069923 ,        -1.4268736907) -- 
(        -.2919069923 ,         1.3029313232) ;
 \node at (         .3292629241 ,         -.0619711837) {    20    }; 
\draw(         .3471338321 ,        -2.4343027495) -- 
(        1.2106435814 ,        -2.4343027495) -- 
(        1.2106435814 ,        -3.2922414682) -- 
(         .3471338321 ,        -3.2922414682) -- 
(         .3471338321 ,        -2.4343027495) ;
 \node at (         .7788887068 ,        -2.8632721088) {     1r}; 
\draw(       -2.3550624466 ,         1.8074128521) -- 
(       -1.3021376555 ,         1.8074128521) -- 
(       -1.3021376555 ,         1.0998919607) -- 
(       -2.3550624466 ,         1.0998919607) -- 
(       -2.3550624466 ,         1.8074128521) ;
 \node at (       -1.8286000510 ,         1.4536524064) {     2r}; 
\draw(       -1.3866483132 ,        -2.4281276269) -- 
(        -.6401302074 ,        -2.4281276269) -- 
(        -.6401302074 ,        -3.5423337550) -- 
(       -1.3866483132 ,        -3.5423337550) -- 
(       -1.3866483132 ,        -2.4281276269) ;
 \node at (       -1.0133892603 ,        -2.9852306910) {     4r}; 
\draw(        -.2728065072 ,        -2.4335318168) -- 
(         .3455778939 ,        -2.4335318168) -- 
(         .3455778939 ,        -3.7817212318) -- 
(        -.2728065072 ,        -3.7817212318) -- 
(        -.2728065072 ,        -2.4335318168) ;
 \node at (         .0363856933 ,        -3.1076265243) {     5r}; 
\draw(       -3.1513382182 ,         -.9257329555) -- 
(       -2.4438173268 ,         -.9257329555) -- 
(       -2.4438173268 ,        -2.4633374123) -- 
(       -3.1513382182 ,        -2.4633374123) -- 
(       -3.1513382182 ,         -.9257329555) ;
 \node at (       -2.7975777725 ,        -1.6945351839) {     7r}; 
\end{tikzpicture}
\captionsetup{font=small}
\caption{Maximising number of rectangles packed, rotation allowed,  container area fraction $\frac{2}{3}$ }
\label{fig3a}
\end{figure}

\begin{figure}[!htb]
\centering
\begin{tikzpicture}[scale=1];
\draw (0,0) circle [radius=        4.0000000000cm];
\draw(         .3721021898 ,         3.3057146847) -- 
(        1.2300409085 ,         3.3057146847) -- 
(        1.2300409085 ,         2.4422049354) -- 
(         .3721021898 ,         2.4422049354) -- 
(         .3721021898 ,         3.3057146847) ;
 \node at (         .8010715491 ,         2.8739598100) {     1    }; 
\draw(        -.3921890695 ,         2.7403462992) -- 
(         .3153318218 ,         2.7403462992) -- 
(         .3153318218 ,         1.6874215081) -- 
(        -.3921890695 ,         1.6874215081) -- 
(        -.3921890695 ,         2.7403462992) ;
 \node at (        -.0384286239 ,         2.2138839036) {     2     }; 
\draw(        1.2682401892 ,         3.2012540707) -- 
(        2.3824463174 ,         3.2012540707) -- 
(        2.3824463174 ,         2.4547359649) -- 
(        1.2682401892 ,         2.4547359649) -- 
(        1.2682401892 ,         3.2012540707) ;
 \node at (        1.8253432533 ,         2.8279950178) {     4     }; 
\draw(        -.5772985576 ,         3.9245494195) -- 
(         .7708908575 ,         3.9245494195) -- 
(         .7708908575 ,         3.3061650184) -- 
(        -.5772985576 ,         3.3061650184) -- 
(        -.5772985576 ,         3.9245494195) ;
 \node at (         .0967961500 ,         3.6153572189) {     5     }; 
\draw(       -2.1980591087 ,         3.3265241580) -- 
(        -.6604546519 ,         3.3265241580) -- 
(        -.6604546519 ,         2.6190032666) -- 
(       -2.1980591087 ,         2.6190032666) -- 
(       -2.1980591087 ,         3.3265241580) ;
 \node at (       -1.4292568803 ,         2.9727637123) {     7     }; 
\draw(       -2.7681367577 ,         2.2836543403) -- 
(       -1.5982203232 ,         2.2836543403) -- 
(       -1.5982203232 ,         1.2697267637) -- 
(       -2.7681367577 ,         1.2697267637) -- 
(       -2.7681367577 ,         2.2836543403) ;
 \node at (       -2.1831785404 ,         1.7766905520) {     8     }; 
\draw(         .4632962526 ,         2.4417576989) -- 
(        2.8978366426 ,         2.4417576989) -- 
(        2.8978366426 ,         1.8735125736) -- 
(         .4632962526 ,         1.8735125736) -- 
(         .4632962526 ,         2.4417576989) ;
 \node at (        1.6805664476 ,         2.1576351362) {     9     }; 
\draw(         .7713044023 ,         1.8730830391) -- 
(        2.3089088591 ,         1.8730830391) -- 
(        2.3089088591 ,          .9538629834) -- 
(         .7713044023 ,          .9538629834) -- 
(         .7713044023 ,         1.8730830391) ;
 \node at (        1.5401066307 ,         1.4134730113) {    10    }; 
\draw(       -2.0432547364 ,         -.1648159916) -- 
(        -.4387979119 ,         -.1648159916) -- 
(        -.4387979119 ,        -1.1843145988) -- 
(       -2.0432547364 ,        -1.1843145988) -- 
(       -2.0432547364 ,         -.1648159916) ;
 \node at (       -1.2410263242 ,         -.6745652952) {    12    }; 
\draw(       -2.3277695664 ,        -1.1961972270) -- 
(        -.2609171987 ,        -1.1961972270) -- 
(        -.2609171987 ,        -2.0151387312) -- 
(       -2.3277695664 ,        -2.0151387312) -- 
(       -2.3277695664 ,        -1.1961972270) ;
 \node at (       -1.2943433826 ,        -1.6056679791) {    13    }; 
\draw(        -.3285471311 ,          .9534326585) -- 
(        2.2396979943 ,          .9534326585) -- 
(        2.2396979943 ,          .1400621849) -- 
(        -.3285471311 ,          .1400621849) -- 
(        -.3285471311 ,          .9534326585) ;
 \node at (         .9555754316 ,          .5467474217) {    14    }; 
\draw(       -3.9208966420 ,          .7135869409) -- 
(       -2.4111473384 ,          .7135869409) -- 
(       -2.4111473384 ,         -.7905913321) -- 
(       -3.9208966420 ,         -.7905913321) -- 
(       -3.9208966420 ,          .7135869409) ;
 \node at (       -3.1660219902 ,         -.0385021956) {    15    }; 
\draw(       -2.6088095747 ,        -2.0261244829) -- 
(         .0541430715 ,        -2.0261244829) -- 
(         .0541430715 ,        -3.0233389676) -- 
(       -2.6088095747 ,        -3.0233389676) -- 
(       -2.6088095747 ,        -2.0261244829) ;
 \node at (       -1.2773332516 ,        -2.5247317252) {    16    }; 
\draw(       -1.4079792775 ,         2.6002441046) -- 
(        -.3940517009 ,         2.6002441046) -- 
(        -.3940517009 ,         -.1462740013) -- 
(       -1.4079792775 ,         -.1462740013) -- 
(       -1.4079792775 ,         2.6002441046) ;
 \node at (        -.9010154892 ,         1.2269850516) {    18    }; 
\draw(       -3.4932893709 ,        -1.0760804855) -- 
(       -2.3512280896 ,        -1.0760804855) -- 
(       -2.3512280896 ,        -1.7668882850) -- 
(       -3.4932893709 ,        -1.7668882850) -- 
(       -3.4932893709 ,        -1.0760804855) ;
 \node at (       -2.9222587303 ,        -1.4214843852) {     3r}; 
\draw(       -2.4109560457 ,         1.0270660988) -- 
(       -1.4081705304 ,         1.0270660988) -- 
(       -1.4081705304 ,          .0354226448) -- 
(       -2.4109560457 ,          .0354226448) -- 
(       -2.4109560457 ,         1.0270660988) ;
 \node at (       -1.9095632880 ,          .5312443718) {     6r    }; 
\draw(       -1.3360063095 ,        -3.0778872641) -- 
(         .9815424370 ,        -3.0778872641) -- 
(         .9815424370 ,        -3.7686950636) -- 
(       -1.3360063095 ,        -3.7686950636) -- 
(       -1.3360063095 ,        -3.0778872641) ;
 \node at (        -.1772319362 ,        -3.4232911638) {    11r    }; 
\draw(        -.3938376836 ,          .1396201047) -- 
(        2.3192542384 ,          .1396201047) -- 
(        2.3192542384 ,         -.8687364413) -- 
(        -.3938376836 ,         -.8687364413) -- 
(        -.3938376836 ,          .1396201047) ;
 \node at (         .9627082774 ,         -.3645581683) {    17r    }; 
\draw(        2.3194623454 ,          .6434938171) -- 
(        3.7567882507 ,          .6434938171) -- 
(        3.7567882507 ,        -1.3732192748) -- 
(        2.3194623454 ,        -1.3732192748) -- 
(        2.3194623454 ,          .6434938171) ;
 \node at (        3.0381252981 ,         -.3648627288) {    19r    }; 
\draw(         .0688138787 ,         -.8702316402) -- 
(        1.8236885305 ,         -.8702316402) -- 
(        1.8236885305 ,        -3.0763597739) -- 
(         .0688138787 ,        -3.0763597739) -- 
(         .0688138787 ,         -.8702316402) ;
 \node at (         .9462512046 ,        -1.9732957070) {    25r    }; 
\end{tikzpicture}
\captionsetup{font=small}
\caption{Maximising total area of rectangles packed,  rotation allowed,   container area fraction $\frac{2}{3}$ }
\label{fig4a}
\end{figure}
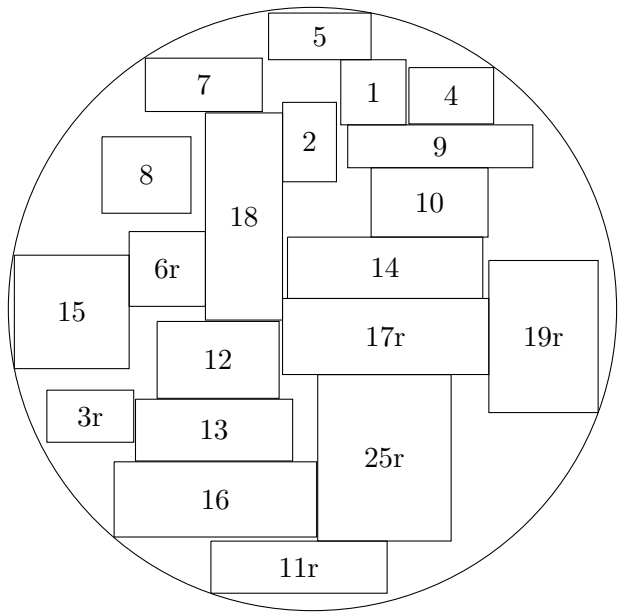

 \clearpage
\newpage
 \pagestyle{empty}
\linespread{1}
\small \normalsize 
\section*{\textbf{Acknowledgments}}
The first author has a grant support from the programme UNAM-DGAPA-PAPIIT-IA106916

\end{document}